\documentclass[12pt]{article}
\usepackage{amssymb, url}
\usepackage{amsmath,amsthm}
\usepackage{mathtools}
\usepackage{dsfont}

\usepackage{verbatim}
\usepackage{graphicx}
\usepackage{epstopdf}
\usepackage{fullpage}
\usepackage[hang,flushmargin]{footmisc}

\usepackage{multicol}
\usepackage{datetime}
\usepackage{bbm}
\usepackage{tikz}
\usepackage{enumitem}
\usepackage{authblk}
\usepackage{color}

\newtheorem{theorem}{Theorem}[]
\newtheorem{corollary}[theorem]{Corollary}
\newtheorem{lemma}{Lemma}
\newtheorem{conjecture}{Conjecture}
\newtheorem{claim}{Claim}[section]

\newtheorem{question}{Question}

\newtheorem*{main1factors}{Theorem~\ref{thm:main-1-factors}}
\newtheorem*{mainmaxdeg}{Theorem~\ref{thm:maxDegreeBound}}
\newtheorem*{mainmiddegThm}{Theorem~\ref{thm:middegreebound}}

\theoremstyle{definition}

\providecommand{\keywords}[1]
{
  \small	
  \textbf{\textit{Keywords---}} #1
}

\begin{document}
\title{On a conjecture that strengthens Kundu's $k$-factor Theorem}
\author{James M.\ Shook$^{3,4}$}

\footnotetext[3]{National Institute of Standards and Technology, Computer Security Division, Gaithersburg, MD; \\ {\tt james.shook@nist.gov}.}
\footnotetext[4]{This article is a U.S. Government work and is not subject to copyright in the USA.}

\maketitle
\begin{abstract}
    Let $\pi=(d_{1},\ldots,d_{n})$ be a non-increasing degree sequence with even $n$. In 1974, Kundu showed that if $\mathcal{D}_{k}(\pi)=(d_{1}-k,\ldots,d_{n}-k)$ is graphic, then some realization of $\pi$ has a $k$-factor. For $r\leq 2$, Busch et al.\ and later Seacrest for $r\leq 4$ showed that if $r\leq k$ and $\mathcal{D}_{k}(\pi)$ is graphic, then there is a realization with a $k$-factor whose edges can be partitioned into a $(k-r)$-factor and $r$ edge-disjoint $1$-factors. We improve this to any $r\leq \min\{\lceil\frac{k+5}{3}\big\rceil,k\}$. In 1978, Brualdi and then Busch et al.\ in 2012, conjectured that $r=k$. The conjecture is still open for $k\geq6$. However, Busch et al.\ showed the conjecture is true when $d_{1}\leq \frac{n}{2}+1$ or $d_{n}\geq \frac{n}{2}+k-2$. We explore this conjecture by first developing new tools that generalize edge-exchanges. With these new tools, we can drop the assumption $\mathcal{D}_{k}(\pi)$ is graphic and show that if $d_{d_{1}-d_{n}+k}\geq d_{1}-d_{n}+k-1,$ then $\pi$ has a realization with $k$ edge-disjoint $1$-factors. From this we confirm the conjecture when $d_{n}\geq \frac{d_{1}+k-1}{2}$ or when $\mathcal{D}_{k}(\pi)$ is graphic and $d_{1}\leq \max \{n/2+d_{n}-k,(n+d_{n})/2\}$. 
\end{abstract}

\section{Introduction}
   
For an undirected graph $G=(V,E)$ with vertex set $V=\{v_{1},\ldots,v_{n}\}$ and edge set $E$, we let $(deg_{G}(v_{1}),\ldots,deg_{G}(v_{n}))$ denote a degree sequence of $G$. We say a sequence $\pi=(d_{1},\ldots, d_{n})$ is graphic if it is the degree sequence of some graph, and call that graph a realization of $\pi$. We let $\mathcal{R}(\pi)$ be the set of realizations of $\pi$, and we let $\pi(G)$ be a degree sequence of a graph $G$ and shorten $\mathcal{R}(\pi(G))$ to $\mathcal{R}(G)$. We say a degree sequence $(d_{1},\ldots,d_{n})$ is non-increasing if $d_{1}\geq \ldots \geq d_{n}$ and positive if $d_{i}\geq 1$ for all $i$. In this paper we assume all degree sequences are non-increasing and only consider graphs and realizations that have no loops or multi-edges. 

In 1974, Kundu \cite{Kundu1974}, followed by Chen \cite{Chen1988} in 1988 with a short proof, gave necessary and sufficient conditions for a degree sequence to have a realization with a spanning near regular subgraph. We call a spanning $k$-regular subgraph a $k$-factor. Since this paper is only concerned with $k$-factors, we present the regular case in Theorem~\ref{kundu}.

\begin{theorem}[Regular case of Kundu's k-factor Theorem \cite{Kundu1974}]\label{kundu}Some realization of a degree sequence $(d_{1},\ldots,d_{n})$ has a $k$-factor if and only if $(d_{1}-k,\ldots,d_{n}-k)$ is graphic.
\end{theorem}

For a sequence $\pi=(d_{1},\ldots,d_{n})$, we let $\mathcal{D}_{k}(\pi)$ denote the sequence $(d_{1}-k,\ldots,d_{n}-k)$ and let $\overline{\pi}=(n-1-d_{n},\ldots,n-1-d_{1})$. Busch, Ferrara, Hartke, Jacobson, Kaul, and West \cite{Busch2012} showed that if both $\pi$ and $\mathcal{D}_{k}(\pi)$ are graphic, then for $r\leq \min\{3,k\}$, there is a realization of $\pi$ with a $k$-factor that has $r$ edge-disjoint $1$-factors. Later, Seacrest \cite{Seacrest2021} improved this to $r\leq \min\{5,k\}$. This naturally leads one to the following question.

\begin{question}\label{que:main}For graphic sequences $\pi$ and $\mathcal{D}_{k}(\pi)$, what is the largest non-negative integer $r\leq k$ such that $\pi$ has a realization with a $k$-factor that has $r$ edge-disjoint $1$-factors?
\end{question}

Besides the recent work of Seacrest and that of Shook \cite{Shook2024} for large $k$, no progress has been made on Question~\ref{que:main} since Busch et al.\ \cite{Busch2012}. We take a fresh look at the problem by building upon some existing strategies by incorporating a generalization of edge-exchanges. A simple edge-exchange (often referred to as a switch or swap) exchanges the edges $uv$ and $xy$ of a graph $G$ with edges $ux$ and $vy$ of the complement of $G$. We have further developed a generalization of the edge-exchange (referred to as multi-switches by Seacrest \cite{Seacrest2021}) that preserves the existence of edge-disjoint regular factors. With these strategies and tools, we show that the number of $1$-factors increases linearly with respect to $k$.

\begin{theorem}\label{thm:main-split}Let $\pi=(d_{1},\ldots,d_{n})$ be a non-increasing positive degree sequence with even $n$. For a positive integers $k\leq d_{n}$, if $\mathcal{D}_{k}(\pi)$ is graphic, then there is some $G\in \mathcal{R}(\pi)$ that has a $k$-factor with at least $\min\{\lceil\frac{k+5}{3}\big\rceil,k\}$ edge-disjoint $1$-factors.
\end{theorem}
Theorem~\ref{thm:main-split} is the first result, that we know of, that shows that the number of $1$-factors grows as a function of $k$. When $k\leq 5$, Theorem~\ref{thm:main-split} shows that $r=k$. Brualdi \cite{Brualdi1978} in 1978 and then later, independently,  Busch et al.\ \cite{Busch2012} conjectured that $r=k$ is true for all $k$. 

\begin{conjecture}[\cite{Brualdi1978} and later in \cite{Busch2012}]\label{conj:KunduExpansion}Some realization of a degree sequence $(d_{1},\ldots,d_{n})$ with even $n$ has $k$ edge-disjoint $1$-factors if and only if $(d_{1}-k,\ldots,d_{n}-k)$ is graphic.
\end{conjecture}

The conjecture is aesthetically pleasing and unresolved, but it is also well supported. For instance, we show later in the paper that there are realizations with at least $\frac{2k+1}{3}$ edge-disjoint $1$-factors. However, we don't know if those $1$-factors can be found within a $k$-factor. The apparent difficulty of the conjecture has motivated mathematicians to look for sufficient conditions. In \cite{Busch2012} the authors showed that if $\mathcal{D}_{k}(\pi)$ is graphic and $d_{n}\geq\frac{n}{2}+k-2$ or $d_{1}\leq \frac{n+2}{2}$, then some realization of $\pi$ has $k$ edge-disjoint $1$-factors. We improve these bounds considerably with our Theorems \ref{thm:main-1-factors}, \ref{thm:maxDegreeBound}, and \ref{thm:middegreebound}. 

\begin{theorem}\label{thm:main-1-factors}Let $\pi=(d_{1},\ldots,d_{n})$ be a non-increasing positive degree sequence with even $n$. For a positive integer $k\leq d_{n}$, if \begin{equation}\label{eq:main}d_{d_{1}-d_{n}+k}\geq d_{1}-d_{n}+k-1,
\end{equation} then there is some $G\in \mathcal{R}(\pi)$ that has $k$ edge-disjoint $1$-factors.
\end{theorem}

Theorem~\ref{thm:main-1-factors} shows us that if $d_{n}\geq \frac{d_{1}+k-1}{2}$, then $\mathcal{D}_{k}(\pi)$ is graphic and some realization of $\pi$ has $k$ edge-disjoint $1$-factors. Based on Theorem~\ref{thm:main-1-factors}, we can prove the following Theorem.
\begin{theorem}\label{thm:maxDegreeBound}Let $\pi=(d_{1},\ldots,d_{n})$ be a non-increasing positive degree sequence with even $n$ such that $\mathcal{D}_{k}(\pi)$ is graphic. If \begin{equation}\label{eq:maxDegreebound}d_{n+1-(d_{1}-d_{n}+k)}\leq n-(d_{1}-d_{n}),\end{equation} then $\pi$ has a realization with $k$ edge-disjoint $1$-factors.
\end{theorem}

Our final main result makes use of the modified Durfee number. For a non-increasing degree sequence $\pi$, the modified Durfee number is defined as $m(\pi)=\max\{i:d_{i}\geq i-1\}$. The Durfee number is $m(\pi)-1$, also known as the strong index \cite{Zverovich1992}, and is closely related to the Durfee square of integer partitions.

\begin{theorem}\label{thm:middegreebound}Let $\pi=(d_{1},\ldots,d_{n})$ be a non-increasing positive degree sequence with even $n$ such that $\mathcal{D}_{k}(\pi)$ is graphic. If $d_{\min\{\frac{n}{2},\,m(\pi)-2\}}> \big\lceil\frac{n+3k-8}{2}\big\rceil$ or $\big\lceil\frac{n+5-k}{2}\big\rceil>d_{\max\{\frac{n}{2}+1,\,n+3-m(\overline{D_{k}(\pi)})\}}$, then $\pi$ has a realization with $k$ edge-disjoint $1$-factors.
\end{theorem}

Later, we will demonstrate how Theorem~\ref{thm:maxDegreeBound} and Theorem~\ref{thm:middegreebound} can be utilized to establish the following lemma, which improves the bound on $d_{1}$.
\begin{lemma}\label{cor:upperd1bounds}Let $\pi=(d_{1},\ldots,d_{n})$ be a non-increasing positive degree sequence with even $n$ such that $\mathcal{D}_{k}(\pi)$ is graphic. If \[d_{1}\leq \max \bigg\{\frac{n}{2}+d_{n}-k,\frac{n+d_{n}}{2}\bigg\},\] then there is a realization of $\pi$ with $k$ edge-disjoint $1$-factors.
\end{lemma}

We have organized the paper as follows. In the next section we discuss past work, various strategies to Question~\ref{que:main} (Section~\ref{strategies}), and our contributions (Section~\ref{constributions}). In Section~\ref{sec:Notation} we present notation and the Gallai-Edmonds Structure Theorem. In Section~\ref{sec:ColorExchanges} we present a generalization of the edge-exchange and prove lemmas that we use extensively in our proofs. Theorem~\ref{thm:main-split} is proved in Section~\ref{sec:main-split}, Theorem~\ref{thm:main-1-factors} is proved in Section~\ref{sec:main-1-factors}, and the proof of Theorem~\ref{thm:middegreebound} can be found in Section~\ref{sec:middegreebound}.

\section{Discussion}\label{discussion}
This section is broken into two subsections. In the first subsection we discuss past work and various approaches towards answering Question~\ref{que:main}. In the second subsection we discuss our contributions (including a few short proofs) and offer a couple of suggestions for future study.
 
\subsection{Past Work and Strategies}\label{strategies}

Finding $1$-factors in $k$-regular graphs is well studied \cite{Plummer2007, Plummer1986, Pulleyblank1996}, and those results are certainly a good place to start when thinking about Question~\ref{que:main}. However, not all $k$-regular graphs can be partitioned into $k$ edge-disjoint $1$-factors. At the extreme, where $t$ is a positive integer, the $2t$-regular graph that is the disjoint union of two complete graphs each with $2t+1$ vertices does not have a $1$-factor. Although, for $k\geq 2\lceil \frac{n}{4}\rceil-1$, the well known $1$-factorization conjecture (See \cite{Chetwynd1985} by Chetwynd and Hilton) implies that every $k$-regular graph can be partitioned into $k$ edge-disjoint $1$-factors. The same $2t$-regular graph we mentioned before shows the lower bound on $k$ is best possible. In a fantastic paper, Csaba et al.\ proved the $1$-factorization conjecture for $n$ sufficiently large.
\begin{theorem}[\cite{Csaba2016}]\label{thm:Csaba}
There exists an $n_{0}\in \mathbb{N}$ such that the following holds. Let
$n,k \in \mathbb{N}$ be such that $n\geq n_{0}$ is even and $k\geq 2\lceil \frac{n}{4}\rceil-1$. Then every $k$-regular graph $G$ on $n$ vertices can be decomposed into $k$ edge-disjoint $1$-factors.
\end{theorem}

A positive resolution of the $1$-factorization conjecture would prove Conjecture~\ref{conj:KunduExpansion} for large $k$, and thanks to Csaba et al.\ we know Conjecture~\ref{conj:KunduExpansion} is true for large $k$ and $n$ sufficiently large. 

If $(d_{1},\ldots, d_{n})$ is a non-increasing degree sequence with large enough $d_{n}$, then the associated realizations are edge dense and should have many edge-disjoint $1$-factors. Hartke and Seacrest \cite{Hartke2012} showed that when $d_{n}\geq \frac{n}{2}+2$ there is a realization with $f(d_{n},n)$ edge-disjoint $1$-factors where \[f(d_{n},n)=\bigg\lfloor\frac{d_{n}-2+\sqrt{n(2d_{n}-n-4)}}{4}\bigg\rfloor.\] Our lower bound on $d_{n}$ from Theorem~\ref{thm:main-1-factors} improves their bound when $d_{n}-f(d_{n},n)\geq d_{1}-d_{n}$. Which is true for the vast majority of possible $d_{1}$ for any given $n$ and $d_{n}$. To see this, we can use the rough lower bound $d_{n}>2f(d_{n},n)$ and $d_{n}\geq n/2+2$ to show any $d_{1}\leq 3n/4+3\leq 3d_{n}/2<2d_{n}-f(d_{n},n)$ will do.

One may think increasing the edge-connectivity of $k$-regular graphs would produce many edge-disjoint $1$-factors. The classic example of this idea is by Berge \cite{Berge1962} and expanded on in \cite{Bollobas1985, Katerinis1993, Ma2023, Mattiolo2022, Plesnik1972, Shiu2008, Thomassen2020}.
\begin{theorem}[\cite{Berge1962}]\label{thm:berge} All even ordered $(k-1)$-edge-connected $k$-regular graphs have a $1$-factor.
\end{theorem}

In \cite{Shook2024} we expanded on a result of Edmonds \cite{Edmonds1964} that gave necessary and sufficient conditions for when a degree sequence has a maximally edge-connected realization. A graph with minimum degree $\delta$ is maximally edge-connected if one has to remove at least $\delta$ edges to disconnect the graph. One of the things we showed is that if $G$ is a simple graph with minimum degree two and has a $1$-factor $F$, then there is a realization of $\mathcal{R}(G)$ that is maximally edge-connected with the subgraph $G-E(F)$. This result can be used to require the realization given by Conjecture~\ref{conj:KunduExpansion} to be maximally edge-connected when $d_{n}\geq 2$. In the same paper we proved a more general result that we combined with Theorem~\ref{thm:berge} to show that if $k$ is large enough, then we can find a realization with a $k$-factor that has at least $r$ edge-disjoint $1$-factors.

\begin{theorem}[\cite{Shook2024}]\label{lem:largek}Let $\pi=(d_{1},\ldots,d_{n})$ be a non-increasing degree sequence with even $n$ such that $\mathcal{D}_{k}(\pi)$ is graphic. If $k\geq \frac{d_{1}}{2}+r-1$ or $k\geq n-1-d_{n}+2(r-1)$, then $\pi$ has a realization with a $k$-factor that has $r$ edge-disjoint $1$-factors.
\end{theorem}

If $k\geq \frac{d_{1}}{2}+r-1$, then since $d_{n}\geq k\geq \frac{d_{1}+2r-2}{2}$, Theorem~\ref{thm:main-1-factors} is stronger than Theorem~\ref{lem:largek} when $k\leq 2r-1$ and Theorem~\ref{thm:Csaba} would be stronger for large $k$ and $n$. However, this strategy of finding highly edge-connected $k$-factors may not be enough to prove Conjecture~\ref{conj:KunduExpansion} since Mattiolo \cite{Mattiolo2022} presented $k$-regular $k$-edge-connected graphs that cannot be partitioned into a $2$-factor and $k-2$ $1$-factors.

In this paper, our general strategy for Question~\ref{que:main} begins by reducing, with great care, the number of realizations we are searching. We then, using one of two possible  methods, study the structure of those realizations with heavy use of counting arguments and the lemmas we developed in Section~\ref{sec:ColorExchanges} about edge-exchanges. The first method (Theorem~\ref{thm:middegreebound}) consists of finding a desirable realization that has a $k$-factor that can be partitioned into $1$-factors and $2$-factors. We then study one of the $2$-factors to show it must have a particular structure that gives us our bound. Our second method follows the lead of \cite{Busch2012} by using the Gallai-Edmonds-Structure Theorem (See Section~\ref{sec:Notation}). We use it in a couple of different ways. In Theorem~\ref{thm:main-1-factors} we find a realization with the most $1$-factors and then apply the Gallai-Edmonds-Structure theorem to the rest of the graph. To prove Theorem~\ref{thm:main-split} we find a realization with a $k$-factor containing the most $1$-factors and then apply the structure theorem to the rest of the $k$-factor. 

\subsection{Our Contributions}\label{constributions}

In this section we discuss our main theorems, present and prove their corollaries, and walk through how they can be combined. Along the way we suggest directions for future study.

Along with requiring the connectivity of a graph $G$ and it's complement $\overline{G}$, Ando et al.\ \cite{Ando1999} showed that bounding the difference of the maximum degree and minimum degree of $G$ would yield a $1$-factor in either $G$ or $\overline{G}$. Ignoring the connectivity requirement we are able to show in Theorem~\ref{thm:main-1-factors} that bounding the difference $d_{1}-d_{n}$ for a non-increasing degree sequence $(d_{1}\ldots, d_{n})$  can tell us if there is a realization with many edge-disjoint $1$-factors.

\begin{main1factors}
     Let $\pi=(d_{1},\ldots,d_{n})$ be a non-increasing positive degree sequence with even $n$. For a positive integer $k\leq d_{n}$, if \begin{equation*}d_{d_{1}-d_{n}+k}\geq d_{1}-d_{n}+k-1,
\end{equation*} then there is some $G\in \mathcal{R}(\pi)$ that has $k$ edge-disjoint $1$-factors.
\end{main1factors}

\begin{figure}
    \centering
    \begin{tikzpicture}
        \draw (1,.5) ellipse (1cm and .5cm);
        \draw (1,2.) ellipse (2cm and .5cm);

        \node[] at (1,.5) {$K_{t+1}$};
        \node[] at (1,2.1) {$I_{t+2}$};

        \draw (-.5,2) -- (.6,.5);
        \draw (0,2) -- (.7,.5);
        \draw (2,2) -- (1.3,.5);
        \draw (2.5,2) -- (1.4,.5);
\end{tikzpicture}
    \caption{$K_{t+1}\ast I_{t+2}$}
    \label{fig:split}
\end{figure}
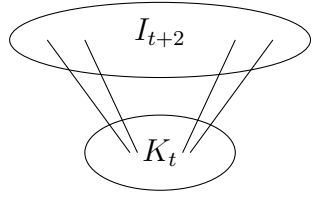

Note that, unlike Conjecture~\ref{conj:KunduExpansion}, in Theorem~\ref{thm:main-1-factors} we did not require $\mathcal{D}_{k}(\pi)$ to be graphic. For $k=1$ and $t\geq 1$, Theorem~\ref{thm:main-1-factors} is best possible since the split graph (See Figure~\ref{fig:split}) joining every vertex of a complete graph $K_{t+1}$ with every vertex of an independent set $I_{t+2}$ has a non-increasing degree sequence $(d_{1},\ldots,d_{n})$ such that $d_{1}=2t+2$, $d_{n}=t+1$, and  $d_{1}-d_{n}=t+1>d_{d_{1}-d_{n}+1}=t$ yet does not have a $1$-factor. However, for $k>1$, we think we can do better. Our motivation for this comes from Corollary~\ref{cor:bestBound}.

\begin{corollary}\label{cor:bestBound}Let $\pi=(d_{1},\ldots,d_{n})$ be a non-increasing positive degree sequence with even $n$. For $k\leq d_{n}$, if \begin{equation}\label{eq:bestbound}d_{d_{1}-d_{n}+1}\geq d_{1}-d_{n}+k-1,\end{equation} then there is some realization of $\pi$ that has a $k$-factor.
\begin{proof}Assume (\ref{eq:bestbound}) is true. The Corollary follows directly from Theorem~\ref{thm:main-1-factors} when $k=1$. Let $t$ be the largest integer such that $\mathcal{D}_{t}(\pi)=(q_{1},\ldots, q_{n})$, where $q_{i}=d_{i}-t$, is graphic. Kundu's $k$-factor theorem implies that $\mathcal{D}_{t+1}(\pi)$ is not graphic, and therefore, no realization of $\mathcal{D}_{t}(\pi)$ has a $1$-factor. This implies $q_{q_{1}-q_{n}+1}<q_{1}-q_{n}$. Since $q_{1}-q_{n}=d_{1}-d_{n}$, we have along with (\ref{eq:bestbound}) that \[d_{d_{1}-d_{n}+1}-t=q_{q_{1}-q_{n}+1}<q_{1}-q_{n}=d_{1}-d_{n}\leq d_{d_{1}-d_{n}+1}-(k-1).\] Which can only be true if $t\geq k$.
\end{proof}
\end{corollary}
Corollary~\ref{cor:bestBound} holds true even when $n$ is odd. However, this result is beyond the scope of our current work, and its non-trivial proof necessitates a different approach. Hence, we defer the proof for a subsequent paper.

Observe that if Conjecture~\ref{conj:KunduExpansion} is true, then Corollary~\ref{cor:bestBound} implies there is a realization with $k$ edge-disjoint $1$-factors. This naturally motivates Conjecture~\ref{con:easierKunduExpansion} as an interesting step towards answering Conjecture~\ref{conj:KunduExpansion}. 

\begin{conjecture}\label{con:easierKunduExpansion}Let $\pi=(d_{1},\ldots,d_{n})$ be a non-increasing positive degree sequence with even $n$. For a positive integer $k\leq d_{n}$, if \begin{equation}\label{eq:mainConj}d_{d_{1}-d_{n}+1}\geq d_{1}-d_{n}+k-1,
\end{equation} then there is some $G\in \mathcal{R}(\pi)$ that has $k$ edge-disjoint $1$-factors.
\end{conjecture}

If we first insist $\mathcal{D}_{k}(\pi)$ is graphic, then we can use Theorem~\ref{thm:main-1-factors} to prove Theorem~\ref{thm:maxDegreeBound}.
\begin{mainmaxdeg}Let $\pi=(d_{1},\ldots,d_{n})$ be a non-increasing positive degree sequence with even $n$ such that $\mathcal{D}_{k}(\pi)$ is graphic. If \begin{equation*}d_{n+1-(d_{1}-d_{n}+k)}\leq n-(d_{1}-d_{n}),\end{equation*} then $\pi$ has a realization with $k$ edge-disjoint $1$-factors.
\begin{proof}Let $q_{i}=d_{i}-k$. We focus on $\mathcal{D}_{k}(\pi)$ and consider its complement $\overline{\mathcal{D}_{k}(\pi)}=(\overline{q}_{1},\dots, \overline{q}_{n})$ where $\overline{q}_{i}=n-1-q_{n+1-i}$.  We have by (\ref{eq:maxDegreebound}) that \[q_{n+1-(q_{1}-q_{n}+k)}=d_{n+1-(d_{1}-d_{n}+k)}-k\leq n-(d_{1}-d_{n})-k=n-1-(q_{1}-q_{n}+k-1).\] From this we can show
\[\overline{q}_{1}-\overline{q}_{n}+k-1=q_{1}-q_{n}+k-1\leq n-1-q_{n+1-(q_{1}-q_{n}+k)}=\overline{q}_{q_{1}-q_{n}+k}=\overline{q}_{\overline{q}_{1}-\overline{q}_{n}+k}.\] Therefore, by Theorem~\ref{thm:main-1-factors}, $\overline{\mathcal{D}_{k}(\pi)}$ has a realization with $k$ edge-disjoint $1$-factors. Thus, those $k$ edge-disjoint $1$-factors can be added to a realization of $\mathcal{D}_{k}(\pi)$ to create a realization of $\pi$ with $k$ edge-disjoint $1$-factors.
\end{proof}
\end{mainmaxdeg}
Note that if Conjecture~\ref{con:easierKunduExpansion} holds, then (\ref{eq:maxDegreebound}) can be improved to $d_{n-(d_{1}-d_{n})}\leq n-(d_{1}-d_{n})$.

\begin{mainmiddegThm}   
Let $\pi=(d_{1},\ldots,d_{n})$ be a non-increasing positive degree sequence with even $n$ such that $\mathcal{D}_{k}(\pi)$ is graphic. If $d_{\min\{\frac{n}{2},m(\pi)-2\}}> \big\lceil\frac{n+3k-8}{2}\big\rceil$ or $\big\lceil\frac{n+5-k}{2}\big\rceil>d_{\max\{\frac{n}{2}+1,n+3-m(\overline{D_{k}(\pi)})\}}$, then $\pi$ has a realization with $k$ edge-disjoint $1$-factors.
\end{mainmiddegThm}

We now combine Theorem~\ref{thm:main-1-factors} and Theorem~\ref{thm:middegreebound}.

\begin{lemma}\label{lem:durfeUpperBound}Let $\pi=(d_{1},\ldots,d_{n})$ be a non-increasing positive degree sequence with even $n$ such that $\mathcal{D}_{k}(\pi)$ is graphic. If \begin{equation}\label{eq:durfeeUpp}d_{1}\leq \frac{m(\pi)+m(\overline{D_{k}(\pi)})}{2}+d_{n}-k,\end{equation} then 
   there is a realization of $\pi$ with $k$ edge-disjoint $1$-factors.
   \begin{proof}We assume (\ref{eq:durfeeUpp}) is true. Let $\overline{\mathcal{D}_{k}(\pi)}=(\overline{q}_{1},\dots, \overline{q}_{n})$ where $\overline{q}_{i}=n-1-d_{n+1-i}+k$ and assume by contradiction that $m(\pi)\leq d_{1}-d_{n}+k-1$ and $m(\overline{D_{k}(\pi)})\leq \overline{q}_{1}-\overline{q}_{n}+k-1=d_{1}-d_{n}+k-1$. Thus, $m(\pi)+m(\overline{D_{k}(\pi)})\leq 2(d_{1}-d_{n}+k-1).$ Solving for $d_{1}$ we have the contradiction \[d_{1}\geq \frac{m(\pi)+m(\overline{D_{k}(\pi)})}{2}+d_{n}-k+1.\] Thus, $m(\pi)> d_{1}-d_{n}+k$ or $m(\overline{D_{k}(\pi)})> \overline{q}_{1}-\overline{q}_{n}+k-1=d_{1}-d_{n}+k$. Applying Theorem~\ref{thm:main-1-factors} to either $\pi$ or $\overline{D_{k}(\pi)}$ we can find a realization of $\pi$ with $k$ edge-disjoint $1$-factors.
   \end{proof}
\end{lemma}

By Theorem~\ref{thm:maxDegreeBound}, $d_{1}\leq \frac{n+d_{n}}{2}$ is enough to show there is a realization of $\pi$ with $k$ edge-disjoint $1$-factors. We can now prove the second bound in 
Lemma~\ref{cor:upperd1bounds} with the help of the following lemma.
 \begin{lemma}[Proved by Li in \cite{Li1975}, but we use the form given by Barrus in \cite{Barrus2022}]\label{lem:durfeeBound}\[
  m(\pi)+m(\overline{\pi}) = 
\begin{cases}
    n+1,& \text{if } d_{m(\pi)}=m(\pi)-1\\
    n,              & \text{otherwise.}
\end{cases}
\] 
\end{lemma}
Since $m(\overline{\pi})\leq m(\overline{D_{k}(\pi)})$,  Lemma~\ref{lem:durfeUpperBound} and Lemma~\ref{lem:durfeeBound} imply
    \[d_{1}\leq\frac{n}{2}+d_{n}-k\leq \frac{m(\pi)+m(\overline{\pi})}{2}+d_{n}-k  \leq\frac{m(\pi)+m(\overline{D_{k}(\pi)})}{2}+d_{n}-k\] is enough to show there is a realization of $\pi$ with $k$ edge-disjoint $1$-factors. Thus, we have proved Lemma~\ref{cor:upperd1bounds}.

If Conjecture~\ref{con:easierKunduExpansion} is true, then we can modify the argument of Lemma~\ref{lem:durfeUpperBound} to show Conjecture~\ref{conj:KunduExpansion} holds for $d_{1}-d_{n}\leq \frac{m(\pi)+m(\overline{D_{k}(\pi)})}{2}-1$. From this we can modify Lemma~\ref{cor:upperd1bounds} to show that $d_{1}-d_{n}\leq \frac{n}{2}-1$ is enough.

If some realization of $(n-1-d_{1},\ldots,n-1-d_{n})$ has a $k'$-factor then we can make use of Petersen's $2$-factor Theorem \cite{Petersen1891} to improve Theorem~\ref{thm:main-split}. Recall that  Petersen showed that any $2q$-regular graph can be partitioned into $q$ $2$-factors.

\begin{theorem}Let $\pi=(d_{1},\ldots,d_{n})$ be a non-increasing positive degree sequence with even $n$. For non-negative integers $r\leq k\leq d_{n}$ and $k'\leq n-1-d_{1}$ such that both $k'$ and $k-r$ are even, if $\mathcal{D}_{k}(\pi)$ and $\pi'=(d_{1}+k',\ldots,d_{n}+k')$ are graphic and \[r\leq \min\bigg\{\bigg\lceil\frac{k+k'+5}{3}\bigg\rceil,k\bigg\},\]  then there is some $G\in \mathcal{R}(\pi)$ that has a $k$-factor with $r$ edge-disjoint $1$-factors.
\begin{proof}Note that $\mathcal{D}_{k+k'}(\pi')=\mathcal{D}_{k}(\pi)$. Since $\pi'$ and $\mathcal{D}_{k+k'}(\pi')$ are graphic and non-increasing we have by Theorem~\ref{thm:main-split} that $\pi'$ has a realization $G$ with a $(k+k')$-factor $F$ with an $r$-factor $F_{0}$ made up of $r$ edge-disjoint $1$-factors. Since $k'+k-r$ is even, we can use Petersen's $2$-factor theorem to split $F-E(F_{0})$ in to $\frac{k'+k-r}{2}$ $2$-factors. We may then select $\frac{k-r}{2}$ of those $2$-factors and add them to $G-E(F)+E(F_{0})$ to construct a realization of $\pi$ with a $k$-factor that has $r$ edge-disjoint $1$-factors.
\end{proof}
\end{theorem}
We can weaken Question~\ref{que:main} by simply asking how many edge-disjoint $1$-factors can a realization of a graphic sequence $\pi$ have if $\mathcal{D}_{k}(\pi)$ is graphic? Seacrest \cite{Seacrest2021} showed that if $\mathcal{D}_{k}(\pi)$ is graphic, then there is a realization of $\pi$ with $\big\lfloor\frac{k}{2}\big\rfloor+2$ edge-disjoint $1$-factors. Seacrest did this by first finding a realization with a $k$-factor that has $r$ edge-disjoint $1$-factors for some $r\equiv k\mod{2}$. In particular, Seacrest used $r=3$ when $k$ is odd and $r=4$ when $k$ is even. Seacrest then took the remaining part of the $k$-factor and applied Petersen's $2$-factor theorem to split it into edge-disjoint $2$-factors. Seacrest then visited each $2$-factor and performed multi-switches, which are defined similarly to the edge-exchanges given in Section~\ref{sec:ColorExchanges} of this paper, to construct at least one additional $1$-factor while leaving existing $1$-factors and $2$-factors intact. This process results in a realization with $\frac{k+r}{2}$ edge-disjoint $1$-factors that may not be part of a $k$-factor. By Theorem~\ref{thm:main-split} there is an $\big\lceil\frac{k+5}{3}\big\rceil-1\leq r\leq \big\lceil\frac{k+5}{3}\big\rceil$ such that $k-r$ is even. This implies that if $\mathcal{D}_{k}(\pi)$ is graphic, then there is a realization of $\pi$ with at least \[\frac{k+r}{2}\geq \frac{k+\big\lceil\frac{k+5}{3}\big\rceil-1}{2}\geq \frac{2k+1}{3}\] edge-disjoint $1$-factors. 

As mentioned before, it was reasonable to pose Conjecture~\ref{conj:KunduExpansion} since one can find a realization with many edge-disjoint $1$-factors when $\mathcal{D}_{k}(\pi)$ is graphic. Perhaps one can adapt the strategies in this paper to find a counter example or affirm the conjecture for $k=6$. For instance, when $k=6$, Theorem~\ref{thm:middegreebound} says we are left to analyze degree sequences for when $\frac{n}{2}+5\geq d_{\frac{n}{2}}\geq d_{\frac{n}{2}+1}\geq\frac{n}{2}$. Perhaps the proof of Theorem~\ref{thm:middegreebound} can be modified to close this gap, or for general $k$, the proof can be modified to find a new bound for $r$ in Question~\ref{que:main}.

\section{Notation and Gallai-Edmonds Decomposition}\label{sec:Notation}

For notation and definitions not defined in this paper, we refer the reader to \cite{Diestel2016}. We let $K_{n}$ denote the complete graph on $n$ vertices. The complement of a graph $G$ is denoted by $\overline{G}$, and $G[X]$ denotes the subgraph induced by $X\subseteq V(G)$. We say a graph is trivial if it has a single vertex. For a graph $G=(V,E)$ and subsets $X$ and $Y$ of $V$, we let $E_{G}(X,Y)$ be the set of all edges in $G$ that have one end in $X$ and the other end in $Y$. We let $e_{G}(X,Y)=|E_{G}(X,Y)|$. For a matching $M$, if $u\in V(M)$, then we let $\overline{u}_{M}$ denote that unique neighbor of $u$ in $M$. Moreover, if $U\subseteq V(M)$ we let $\overline{U}_{M}=\{\overline{u}_{M}\| \forall u\in U\}$. 

Let $G$ be a graph. The number of odd components in $G$ is denoted by $o(G)$. For $S\subseteq V(G)$, we let $\text{def}_{G}(S)=o(G-S)-|S|$, and define the deficiency of $G$ by $\text{def}(G)=\max_{S\subseteq V(G)}\text{def}_{G}(S)$. The Berge-Tutte Formula \cite{Berge1962} says that if $G$ has $n$ vertices, then the maximum size of a matching in $G$ is $\frac{1}{2}(n-\text{def}(G))$. If every subgraph obtained by deleting one vertex from $G$ has a $1$-factor, then we say $G$ is factor-critical. If a matching in $G$ covers all but one vertex, then we say the matching is near-perfect.

In a graph $G$, the Gallai-Edmonds Decomposition of $G$ is a partition of $V(G)$ into three sets $A$, $C$, and $D$ such that $D=V(G)-B$ where $B$ is the set of vertices that are in every maximum matching of $G$ and $B=A\cup C$ where $A$ is the set of vertices of $B$ with at least one neighbor in $D$. The Gallai-Edmonds Structure Theorem (See \cite{West2011} for a short proof and history.) is an important tool for our work and we present it below before beginning our proofs.

\begin{theorem}[Gallai-Edmonds Structure Theorem]\label{thm:GallieEdmonds} Let $A$, $C$, and $D$ be the sets in the Gallai-Edmonds Decomposition of a graph $G$. Let $G_{1}, \ldots, G_{k}$ be the components of $G[D]$. If $M$ is a maximum matching in $G$, then the following properties hold.
\begin{enumerate}[label=(\Roman*),ref=(\Roman*)]
    \item \label{GEST:1} $M$ covers $C$ and matches $A$ into distinct components of $G[D]$.
    \item \label{GEST:2}Each $G_{i}$ is factor-critical, and $M$ restricts to a near-perfect matching on $G_{i}$.
    \item \label{GEST:3}If $\emptyset\neq S\subseteq A$, then $N_{G}(S)$ has a vertex in at least $|S|+1$ of $G_{1},\dots,G_{k}$.
    \item \label{GEST:4}$\text{def}(A)=\text{def}(G)=k-|A|$.
\end{enumerate}
\end{theorem}

\section{Generalized Edge-Exchanges}\label{sec:ColorExchanges}
\begin{figure}[h]
    \centering
    \begin{tikzpicture}
\draw (1,3) -- (2,4);
\draw (3,3) -- (2,2);
\draw[dashed, color=red] (1,3) -- (2,2);
\draw[dashed, color=red] (3,3) -- (2,4);
\draw[fill=black] (1,3) circle (.1cm);
\draw[fill=black] (3,3) circle (.1cm);
\draw[fill=black] (2,2) circle (.1cm);
\draw[fill=black] (2,4) circle (.1cm);
\node[] at (.5,3) {$x_{0}$};
\node[] at (3.5,3) {$x_{1}$};
\node[] at (2.5,4) {$v$};
\node[] at (2.5,2) {$u$};


\draw (10,4) -- (8,3);
\draw (10,2) -- (12,3);
\draw[dashdotted, color=blue] (10,2) -- (8,3);
\draw[dashdotted, color=blue] (10,4) -- (9,3);
\draw[dotted, color=black] (10,2) -- (9,3);
\draw[dotted, color=black] (10,4) -- (11,3);
\draw[dashed, color=red] (10,4) -- (12,3);
\draw[dashed, color=red] (10,2) -- (11,3);

\draw[fill=black] (10,4) circle (.1cm);
\draw[fill=black] (10,2) circle (.1cm);
\draw[fill=black] (9,3) circle (.1cm);
\draw[fill=black] (8,3) circle (.1cm);
\draw[fill=black] (11,3) circle (.1cm);
\draw[fill=black] (12,3) circle (.1cm);

\node[] at (10.4,4) {$v$};
\node[] at (10.4,2) {$u$};
\node[] at (7.6,3) {$x_{0}$};
\node[] at (9.4,3) {$x_{1}$};
\node[] at (10.6,3) {$x_{2}$};
\node[] at (12.4,3) {$x_{3}$};

\end{tikzpicture}
    \caption{Edge-exchanges with length 2 and length 4, respectively.}
    \label{fig:Exchanges}
\end{figure}
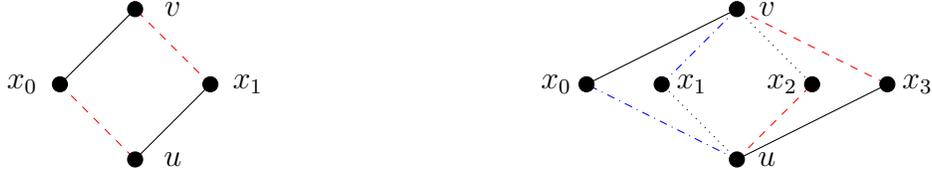

An edge-exchange (see the left side of Figure~\ref{fig:Exchanges}), first discovered by Petersen in his famous 1891 paper \cite{Petersen1891}, consists of exchanging two edges $vx_{0}$ and $x_{1}u$ of $G$ with two edges $x_{0}u$ and $vx_{1}$ of $\overline{G}$.  Edge-exchanges, also commonly called a switch or a swap, are a fundamental operation when studying degree sequences. Fulkerson et al.\ \cite{Fulkerson1965} showed that through a sequence of edge-exchanges one can transform one realization to any other realization. However, we will need a more general form of edge-exchanges and we present it from the perspective of an edge coloring of $K_{n}$.

Consider an edge coloring of $K_{n}$ with natural numbers $\{1,\ldots,t\}$. We let  $H_{1},\ldots, H_{t}$ denote the subgraphs of $K_{n}$ where $H_{j}$ is formed by all edges colored $j$. Given distinct vertices $v$ and $u$, we say the colors of edges $vx_{0}$ an $x_{0}u$ can be exchanged if there exists a natural number $l$ and a list of $2l$ distinct edges 
\[(vx_{0},x_{0}u,vx_{1},x_{1}u,\ldots,vx_{l-1},x_{l-1}u)\] such that $x_{i}u$ and $vx_{i+1}$ have the same color for all $i$ modulo $l$. Indeed if we exchange the colors of $vx_{i}$ and $x_{i}u$ for all $i$ modulo $l$, then we would create another edge coloring of $K_{n}$ with color classes $H'_{1},\ldots, H'_{t}$ such that $H'_{j}\in \mathcal{R}(H_{j})$. The right hand side of Figure~\ref{fig:Exchanges} shows, as an example, the exchange $(vx_{0},x_{0}u,vx_{1},x_{1}u,vx_{2},x_{2}u,vx_{3},x_{3}u)$. We often just say two edges can be exchanged when it is clear we mean exchanging their colors. If $x_{j}u$ and $vx_{0}$ are not the same color for $j\leq l-1$, then we call the list a near exchange with length $l$. For a near exchange or exchange $L=(vx_{0},x_{0}u,vx_{1},x_{1}u,\ldots,vx_{l-1},x_{l-1}u)$, we let $\mathcal{X}(L)=\{x_{0},\ldots,x_{l-1}\}$.

The rest of this section focuses on exchanging edges where the first edge is in $H_{1}$. However, in later sections we will want to consider exchanges that start with an edge of $H_{2}$ so it is important to point out that all of the results in this section still hold when every occurrence of $H_{1}$ or $H_{2}$ are swapped with each other. 

Let $L$ be a $vx_{0}$ and $x_{0}u$ exchange. If for any $H_{i}$ there is at most one $x_{j}\in \mathcal{X}(L)$ such that $vx_{j}\in E(H_{i})$, then we call $L$ simplified.

\begin{lemma}\label{lem:reduce}Let $H_{1},\ldots, H_{t}$ be the subgraphs formed by coloring every edge of $K_{n}$ with some integer in $\{1,\dots,t\}$ such that $H_{j}$ is a spanning regular graph for $j\geq 3$. For edges $vx_{0}\in E(H_{1})$ and $x_{0}u\notin E(H_{1})$, if $L$ is an exchange for $vx_{0}$ and $x_{0}u$, then there exists a simplified $vx_{0}$ and $x_{0}u$ exchange $L'$ with $\mathcal{X}(L')\subseteq \mathcal{X}(L)$.
\begin{proof}Let $L=(vx_{0},x_{0}u,vx_{1},x_{1}u,\ldots,vx_{l-1},x_{l-1}u)$ be the shortest counter example. Thus, there is an $x_{j}$ and an $x_{t}$ in $\mathcal{X}(L)$ with $j< t$ such that both $vx_{j}$ and $vx_{t}$ are in $E(H_{i})$ for some $i$. We have a contradiction since we can create the shorter exchange \[L'=(vx_{0},x_{0}u,\ldots, vx_{j-1},x_{j-1}u, vx_{t}, x_{t}u,\ldots, vx_{l-1},x_{l-1}u)\] with $\mathcal{X}(L')\subseteq \mathcal{X}(L)$.
\end{proof}
\end{lemma}

\begin{lemma}\label{lem:nearExchange}Let $H_{1},\ldots, H_{t}$ be the subgraphs formed by coloring every edge of $K_{n}$ with some integer in $\{1,\dots,t\}$ such that $H_{j}$ is a spanning regular graph for $j\geq 3$. For edges $vx_{0}\in E(H_{1})$ and $x_{0}u\notin E(H_{1})$, if $vx_{0}$ and $x_{0}u$ cannot be exchanged, then a longest near edge-exchange using $vx_{0}$ and $x_{0}u$ ends with an edge of $H_{2}$.
\begin{proof}Let $L=(vx_{0},x_{0}u,vx_{1},x_{1}u,\ldots,vx_{l-1},x_{l-1}u)$ be a longest near edge-exchange. If there is some $x_{j}u\in E(H_{1})$, then $(vx_{0},x_{0}u,\ldots, vx_{j},x_{j}u)$ would be an edge-exchange since $vx_{0}\in E(H_{1})$.  Suppose $x_{l-1}u\in E(H_{j})$ for some $j\geq 3$. Since $H_{j}$ is regular and $u$ and $v$ are incident to the same number of edges of $H_{j}$ in $\{vx_{0},x_{0}u,\ldots, vx_{l-1},x_{l-1}u\}$ there must be an $x_{l}\in N_{H_{j}}(v)-\mathcal{X}(L)$. However, $(vx_{0},x_{0}u,\ldots, vx_{l-1},x_{l-1}u,vx_{l},x_{l}u)$ would be a longer near edge-exchange contradicting our choice of $l$.
\end{proof}
\end{lemma}

\begin{lemma}\label{lem:guaranteExchange}Let $H_{1},\ldots, H_{t}$ be the subgraphs formed by coloring every edge of $K_{n}$ with some integer in $\{1,\dots,t\}$ such that $H_{j}$ is a spanning regular graph for $j\geq 3$. For edges $vx_{0}\in E(H_{1})$ and $x_{0}u\notin E(H_{1})$, if there is a $y\in N_{H_{2}}(v)\cap N_{H_{1}}(u)$,  or if a $y\in N_{H_{2}}(v)- N_{H_{2}}(u)$ and a $y'\in N_{H_{1}}(u)-N_{H_{1}}(v)$ such that $yu$ and $vy'$ have the same color, then $vx_{0}$ and $x_{0}u$ can be exchanged.
\begin{proof}Suppose $vx_{0}$ and $x_{0}u$ cannot be exchanged. Let $(vx_{0},x_{0}u,\ldots,vx_{l-1},x_{l-1}u)$ be a longest near exchange. By Lemma~\ref{lem:nearExchange} $x_{l-1}u\in E(H_{2})$, and thus, there is a smallest $j$ such that $x_{j}u\in E(H_{2})$. However, we have a contradiction since $(vx_{0},x_{0}u,\ldots,vx_{j},x_{j}u,vy,yu)$ is an exchange when $y\in N_{H_{2}}(v)\cap N_{H_{1}}(u)$, and otherwise,  $(vx_{0},x_{0}u,\ldots,vx_{j},x_{j}u,vy,yu,vy',y'u)$ is an exchange.
\end{proof}
\end{lemma}

\begin{lemma}\label{lem:largedegree}Let $H_{1},\ldots, H_{t}$ be the subgraphs formed by coloring every edge of $K_{n}$ with some integer in $\{1,\dots,t\}$ such that $H_{j}$ is a spanning regular graph for $j\geq 3$. For vertices $u$ and $v$, let $X=\{x^{(1)}_{0},\ldots,x^{(|X|)}_{0}\}$ where $X\subseteq N_{H_{1}}(v)-N_{H_{1}}(u)$.  If \begin{equation}\label{eq:largedegree}deg_{H_{1}}(u)\geq deg_{H_{1}}(v)-|N_{H_{2}}(u)\cap N_{H_{1}}(v)|+|X\cap N_{H_{2}}(u)|,
\end{equation} then there exists a set $\mathcal{L}=\{L^{(1)},\ldots,L^{(|X|)}\}$  such that $L^{(j)}\in \mathcal{L}$ is a $vx^{(j)}_{0}$ and $x^{(j)}_{0}u$ exchange and \[\mathcal{X}(L^{(j)})\cap \mathcal{X}(L^{(i)})=\emptyset\] for $j\neq i$.
\begin{proof}Trivially, $\mathcal{X}((vx^{(j)}_{0},x^{(j)}_{0}u))=\{x_{0}^{(j)}\}$. Thus, there exists a set $\{L^{(1)},\ldots,L^{(|X|)}\}$ and an $1\leq f \leq |X|$, where \[L^{(j)}=(vx^{(j)}_{0},x^{(j)}_{0}u,vx^{(j)}_{1},x^{(j)}_{1}u,\ldots,vx^{(j)}_{l^{(j)}-1},x^{(j)}_{l^{(j)}-1}u)\] is a simplified exchange for $j < f$ and a near exchange for $j\geq f$ with $\mathcal{X}(L^{(j)})\cap \mathcal{X}(L^{(i)})=\emptyset$ for all $i\neq j$, such that \begin{equation}\label{eq:existsExchange}\sum_{i=1}^{|X|}|\mathcal{X}(L^{(i)})|\end{equation} is maximized.

Let $Y=\bigcup_{i=1}^{|X|}\mathcal{X}(L^{(i)})$. Suppose there exists an $s\geq f$ such that $x^{(s)}_{l^{(s)}-1}u\in E(H_{i})$ for some $i\geq 3$. Since $x^{(j)}_{t}u$ and $vx^{(j)}_{t+1}$  are the same color for all $j$ and $t\leq l^{(j)}-2$, we know for each $L^{(j)}$ that $u$ is incident with at least as many edges of $H_{i}$ than $v$. Since $H_{i}$ is regular and $u$ is incident with one more edge of $H_{i}$ than $v$ in $L^{(s)}$, there must be an $x\in N_{H_{i}}(v)-Y$. However, $(vx^{(s)}_{0},x^{(s)}_{0}u,\ldots, vx^{(s)}_{l^{(s)}-1},x^{(s)}_{l^{(s)}-1}u,vx,xu)$ would be a longer near edge-exchange. This contradicts (\ref{eq:existsExchange}) being maximized. Thus, $x^{(j)}_{l^{(j)}-1}u\in E(H_{2})$ for every $j\geq f$. 

We can rewrite (\ref{eq:largedegree}) so that \[deg_{H_{2}}(v)\geq deg_{H_{2}}(u)-|N_{H_{2}}(u)\cap N_{H_{1}}(v)|+|X\cap N_{H_{2}}(u)|=|N_{H_{2}}(u)-N_{H_{1}}(v)|+|X\cap N_{H_{2}}(u)|.\]
Since each $L^{(j)}$ is a simplified exchange or a near exchange, there is at most one $x^{(j)}_{i}\in \mathcal{X}(L^{(j)})$ such that $vx^{(j)}_{i}\in E(H_{1})$.  Thus, $Y\cap N_{H_{1}}(v)=X$, and  therefore, $(Y\cap N_{H_{2}}(u)) \cap N_{H_{1}}(v)=X\cap N_{H_{2}}(u)$ and $(Y\cap N_{H_{2}}(u))-N_{H_{1}}(v)\subseteq N_{H_{2}}(u)-N_{H_{1}}(v)$. Moreover, \[|\mathcal{X}(L^{(j)})\cap N_{H_{2}}(u)|= |\mathcal{X}(L^{(j)})\cap N_{H_{2}}(v)|\] for $j<f$ and \[|\mathcal{X}(L^{(j)})\cap N_{H_{2}}(u)|= |\mathcal{X}(L^{(j)})\cap N_{H_{2}}(v)|+1\] for $j\geq f$. Thus, $|Y\cap N_{H_{2}}(u)|>|Y\cap N_{H_{2}}(v)|$. We may now deduce that there is an $x\in N_{H_{2}}(v)-Y$ since \[deg_{H_{2}}(v)\geq |N_{H_{2}}(u)-N_{H_{1}}(v)|+|X\cap N_{H_{2}}(u)|\geq |Y\cap N_{H_{2}}(u)|>|Y\cap N_{H_{2}}(v)|.\] However, we have a contradiction to (\ref{eq:existsExchange}) since $(vx^{(f)}_{0},x^{(f)}_{0}u,\ldots, vx^{(f)}_{l^{(f)}-1},x^{(f)}_{l^{(f)}-1}u,vx,xu)$ would be a longer near exchange than $L^{(f)}$.
\end{proof}
\end{lemma}

\begin{lemma}\label{lem:existsAnExch}Let $H_{1},\ldots, H_{t}$ be the subgraphs formed by coloring every edge of $K_{n}$ with some integer in $\{1,\dots,t\}$ such that $H_{j}$ is a spanning regular graph for $j\geq 3$. For vertices $u$ and $v$, let $X=\{x^{(1)}_{0},\ldots,x^{(|X|)}_{0}\}$ where $X\subseteq N_{H_{1}}(v)-N_{H_{1}}(u)$. If \begin{equation}\label{eq:existsAnExch}|X-N_{H_{2}}(u)|+|N_{H_{2}}(v)-N_{H_{2}}(u)|>|N_{H_{2}}(u)-N_{H_{1}}(v)-N_{H_{2}}(v)|,
    \end{equation}then there exists an $x^{(j)}_{0}\in X$ such that $vx^{(j)}_{0}$ and $x^{(j)}_{0}u$ can be color exchanged.
\begin{proof}Suppose the lemma is not true. Trivially, $\mathcal{X}((vx^{(j)}_{0},x^{(j)}_{0}u))=\{x_{0}^{(j)}\}$. Thus, there exists a set $\{L^{(1)},\ldots,L^{(|X|)}\}$, where \[L^{(j)}=(vx^{(j)}_{0},x^{(j)}_{0}u,vx^{(j)}_{1},x^{(j)}_{1}u,\ldots,vx^{(j)}_{l^{(j)}-1},x^{(j)}_{l^{(j)}-1}u)\] is a near exchange for $1\leq j\leq |X|$ with $\mathcal{X}(L^{(j)})\cap \mathcal{X}(L^{(i)})=\emptyset$ for all $i\neq j$, such that \begin{equation}\label{eq:existsExchange2}\sum_{i=1}^{|X|}|\mathcal{X}(L^{(i)})|\end{equation} is maximized.
Let $Y=\bigcup_{i=1}^{|X|}\mathcal{X}(L^{(i)})$. 
Suppose there exists an $s$ such that $x^{(s)}_{l^{(s)}-1}u\in E(H_{i})$ for some $i\geq 3$. Since $x^{(j)}_{t}u$ and $vx^{(j)}_{t+1}$  are the same color for all $j$ and $t$, we know for each $L^{(j)}$ that $u$ is incident with at least as many edges of $H_{i}$ than $v$. Since $H_{i}$ is regular and $u$ is incident with one more edge of $H_{i}$ than $v$ in $L^{(s)}$, there must be an $x\in N_{H_{i}}(v)-Y$. However, $(vx^{(s)}_{0},x^{(s)}_{0}u,\ldots, vx^{(s)}_{l^{(s)}-1},x^{(s)}_{l^{(s)}-1}u,vx,xu)$ would be a longer near edge-exchange. This contradicts (\ref{eq:existsExchange2}) being maximized. Thus, $x^{(j)}_{l^{(j)}-1}u\in E(H_{2})$ for every $j$.

If there is an $x\in N_{H_{2}}(v)-N_{H_{2}}(u)$ not in $Y$, then \[(vx^{(1)}_{0},x^{(1)}_{0}u,vx^{(1)}_{1},x^{(1)}_{1}u,\ldots,vx^{(1)}_{l^{(1)}-1},x^{(1)}_{l^{(1)}-1}u,vx, xu)\] is a longer near exchange that contradicts (\ref{eq:existsExchange2}). Thus, $N_{H_{2}}(v)-N_{H_{2}}(u)\subseteq Y$, and therefore, $(X-N_{H_{2}}(u))\cup (N_{H_{2}}(v)-N_{H_{2}}(u))\subseteq Y$.

By (\ref{eq:existsAnExch}) $(X-N_{H_{2}}(u))\cup (N_{H_{2}}(v)-N_{H_{2}}(u))$ is not empty. Consider an $x^{(s)}_{j}\in (X-N_{H_{2}}(u))\cup (N_{H_{2}}(v)-N_{H_{2}}(u))$. Since $x^{(s)}_{l^{(s)}-1}\in N_{H_{2}}(u)$, there is a smallest $f\geq 1$ such that $x^{(s)}_{j+f}\in N_{H_{2}}(u)$. Choose an $x^{(s)}_{j+i}$ for some $1\leq i \leq f$.  By the minimality of $f$ and the fact that $x^{(s)}_{j}\notin N_{H_{2}}(u)$, we know that $x^{(s)}_{j+i-1}\notin N_{H_{2}}(u)$, and therefore, $x^{(s)}_{j+i}\notin N_{H_{2}}(v)-N_{H_{2}}(u)$. Moreover, the definition of a near exchange says $x^{(s)}_{t}\notin N_{H_{1}}(v)$ for all $0<t\leq l^{(s)}-1$. Thus, both $x^{(s)}_{j+f}$ and $x^{(s)}_{j+i}$ are not in $N_{H_{1}}(v)$. Therefore, $x^{(s)}_{j+f}\in N_{H_{2}}(u)-N_{H_{1}}(v)-N_{H_{2}}(v)$, and for $1\leq i \leq f$, $x^{(s)}_{j+i} \notin (X-N_{H_{2}}(u))\cup (N_{H_{2}}(v)-N_{H_{2}}(u))$. Thus, for every vertex in $(X-N_{H_{2}}(u))\cup (N_{H_{2}}(v)-N_{H_{2}}(u))$ there is a vertex in  $N_{H_{2}}(u)-N_{H_{1}}(v)-N_{H_{2}}(v)$. However, this implies the contradiction \[|N_{H_{2}}(u)-N_{H_{1}}(v)-N_{H_{2}}(v)|\geq |(X-N_{H_{2}}(u))\cup (N_{H_{2}}(v)-N_{H_{2}}(u))|.\] Thus, we have proved our Lemma. \qedhere
\end{proof}
\end{lemma}

\section{Proof of Theorem~\lowercase{\ref{thm:main-1-factors}}}\label{sec:main-1-factors}

\begin{proof}We assume every $G\in \mathcal{R}(\pi)$ has vertex set $V=\{v_{1},\ldots,v_{n}\}$ such that $deg_{G}(v_{i})=d_{i}$. Let $r\leq k$ be the largest integer such that there is a realization of $\pi$ with $r$ edge-disjoint $1$-factors. By contradiction we assume $r\leq k-1$. 
Let $\mathcal{G}$ be the set of tuples of the form $(G,F,t)$ where $G\in \mathcal{R}(\pi)$, $F$ is an $r$-factor of $G$ whose edges can be partitioned into $r$ $1$-factors, and $v_{t}\notin V(M)$ for some maximum matching $M$ of $G-E(F)$.

\begin{enumerate}[label=(C\arabic*),ref=(C\arabic*)]
    \item\label{highDegreeMatching_choice:1} We choose a $(G,F,t)\in \mathcal{G}$ such that $\text{def}(G-E(F))$ is minimized, and
    \item\label{highDegreeMatching_choice:2} subject to $\ref{highDegreeMatching_choice:1}$, we minimize $t$. 
\end{enumerate}

Let $M_{1},\ldots, M_{r}$ be a partition of $E(F)$ into $r$ edge-disjoint $1$-factors. Furthermore, we let $M$ be a maximum matching of $G-E(F)$ that misses $v_{t}$. Let $Q=\{v_{1},\ldots,v_{t-1}\}$.

Let $H_{1}=G-E(F)$ and $H_{2}=\overline{G}$, and note that $H_{1},H_{2}, H_{3},\ldots, H_{r+2}$ where $H_{i}=M_{i-2}$ for $i\geq 3$ represent a coloring of the edges of $K_{n}$. Thus, any edge-exchange involving any $H_{i}$ corresponds to a $(G',F',t')\in \mathcal{G}$. 

Let $A$, $C$, and $D$ be a Gallai-Edmonds Decomposition of $H_{1}$. We know that $D$ is not empty since our assumption is that $H_{1}$ does not have a matching. We let $D'\subseteq D-\{v_{t}\}$ be the largest set such that for every $u\in D'$ there is a matching $M_{u}$ in $H_{1}$ that misses both $u$ and $v_{t}$ such that $E(M_{u}-V(D'))=E(M-V(D'))$. Since $M$ misses $v_{t}$ and some other vertex, $D'$ is not empty. Let $A'\subseteq A$ be all vertices in $A$ that are adjacent in $H_{1}$ to a vertex in $D'$. 

\begin{claim}\label{cl:inA}$N_{H_{1}}(u)\subseteq A'\subseteq Q$ for all $u\in D'$.
\begin{proof}By our selection of $v_{t}$ in \ref{highDegreeMatching_choice:2}, we know that $Q\cap D'=\emptyset$. Suppose there is a $v\in N_{H_{1}}(u)$ not in $Q$. Consider a maximum matching $M_{u}$ of $H_{1}$ that misses both $u$ and $v_{t}$. By Lemma~\ref{lem:guaranteExchange} we may exchange $v_{t}u$ and $uv$ to find a $(G',F',t')\in \mathcal{G}$. However, we have a contradiction since $M_{u}+\{v_{t}u\}$ is a maximum matching in $G'-F'$ that violates $\ref{highDegreeMatching_choice:1}$. Thus, $v\in Q$, and therefore, $N_{H_{1}}(u)\subseteq A'\subseteq Q$ for all $u\in D'$. 
\end{proof}
\end{claim}

\begin{claim}\label{cl:A'matchesD'}$\overline{w}_{M}\in D'$ for every $w\in A'$
\begin{proof}
Suppose there is a $w\in A'$ such that $\overline{w}_{M}\in D-D'$. By the definition of $D'$, for any $u\in N_{H_{1}}(w)\cap D'$ there is a $M_{u}$ that misses $u$ and $v_{t}$ such that $E(M_{u}-V(D'))=E(M-V(D'))$. However, since $M'=M_{u}-\{w\overline{w}_{M}\}+\{wu\}$ is a maximum matching that misses $v_{t}$ with $E(M'-V(D'\cup \{w_{M}\}))=E(M-V(D'\cup \{w_{M}\}))$, we see that $D'\cup \{w_{M}\}$ is a larger set than $D'$.
\end{proof}
\end{claim}

Claim~\ref{cl:inA} implies every component of $H_{1}[D']$ is a single vertex. By Claim~\ref{cl:A'matchesD'} and \ref{GEST:3} of Theorem~\ref{thm:GallieEdmonds}, we deduce that $|D'|>|A'|$. Therefore, \[e_{H_{1}}(D',A')\geq |D'|(d_{n}-r)>|A'|(d_{n}-r),\] and by the pigeon hole principle there is some vertex $s\in A'$ adjacent in $H_{1}$ to at least $d_{n}-r+1$ vertices in $D'$.

\begin{claim}\label{cl:QCompletetoA}$Q$ is complete in $G$ to $A'\cup N_{H_{1}}(v_{t})$.
\begin{proof}Suppose there is a $w\in Q$ and $v\in A'\cup N_{H_{1}}(v_{t})$ that are not adjacent in $G$, and let $u\in N_{H_{1}}(v)\cap D'$. By the definition of $D$, there is a maximum matching $M'$ of $H_{1}$ that misses $u$. By \ref{GEST:2} of Theorem~\ref{thm:GallieEdmonds}, $u$ and $\overline{w}_{M'}$ do not occupy the same component of $H_{1}[D]$, and therefore, $u\overline{w}_{M'}\notin E(H_{1})$. By Lemma~\ref{lem:reduce} and Lemma~\ref{lem:guaranteExchange}, there exists a simplified $\overline{w}_{M'}w$ and $u\overline{w}_{M'}$ exchange that when exchanged creates another $(G',F',t')\in \mathcal{G}$ with the matching  $M'-\{\overline{w}_{M'}w\}+\{\overline{w}_{M'}u\}$ of $G'-E(F')$ that violates \ref{highDegreeMatching_choice:2}.
\end{proof}
\end{claim}

Since $s\in A'\subseteq Q$, Claim~\ref{cl:QCompletetoA} states that $s$ is adjacent in $G$ to every vertex in $Q-\{s\}\cup N_{H_{1}}(v_{t})$. Thus, \begin{align}
 |Q-\{s\}|&\leq deg_{G}(s)-|N_{H_{1}}(s)-Q|-|N_{F}(s)-Q|\notag\\
 &\leq d_{1}-|N_{H_{1}}(s)\cap D'|-|N_{H_{1}}(s)\cap \{v_{t}\}|-|N_{F}(s)-Q|\notag\\
 &\leq d_{1}-(d_{n}-r+1)-|N_{H_{1}}(s)\cap \{v_{t}\}|-|N_{F}(s)-Q|\notag\\
&=d_{1}-d_{n}+r-1-|N_{H_{1}}(s)\cap \{v_{t}\}|-|N_{F}(s)-Q|.   \notag
\end{align}
Since $s\in Q$, we deduce that $|Q|\leq d_{1}-d_{n}+r-|N_{H_{1}}(s)\cap \{v_{t}\}|-|N_{F}(s)-Q|$. However, since $|Q|=t-1$, we deduce from (\ref{eq:main}) that \[d_{t}=d_{|Q|+1}\geq d_{d_{1}-d_{n}+r+1}\geq d_{d_{1}-d_{n}+k}\geq d_{1}-d_{n}+k-1\geq d_{1}-d_{n}+r.\]

Suppose $N_{H_{1}}(s)\cap D'\subseteq N_{H_{2}}(v_{t})$. Let $X=\{x\}$ for some $x\in N_{H_{1}}(s)\cap D'$. We let $M_{x}$ be a maximum matching in $H_{1}$ that misses both $x$ and $v_{t}$. Since $deg_{H_{1}}(v_{t})\geq d_{1}-d_{n}$, $|N_{H_{2}}(v_{t})\cap N_{H_{1}}(s)|\geq d_{n}-r+1$, and $|X\cap N_{H_{2}}(v_{t})|=1$, we deduce that \[deg_{H_{1}}(s)-|N_{H_{2}}(v_{t})\cap N_{H_{1}}(s)|+|X\cap N_{H_{2}}(v_{t})|\leq d_{1}-r-(d_{n}-r+1)+1=d_{1}-d_{n}\leq deg_{H_{1}}(v_{t}).\] Therefore, by Lemma~\ref{lem:largedegree} and Lemma~\ref{lem:reduce}, there exist a simplified $v_{t}x$ and $xs$ exchange that when exchanged creates another $(G',F',t')\in \mathcal{G}$ such that $M_{x}$ is a matching of $G'-E(F')$. However, this violates \ref{highDegreeMatching_choice:1} since $M_{x}+\{v_{t}x\}$ is a larger matching of $G'-E(F')$. Thus, $N_{H_{1}}(s)\cap D'\cap N_{F}(v_{t})\neq \emptyset$.

Let $X=N_{H_{1}}(s)\cap D'\cap N_{F}(v_{t})$, and observe that $N_{F}(s)\cap N_{H_{2}}(v_{t})=N_{H_{2}}(v_{t})-N_{H_{1}}(s)-N_{H_{2}}(s)$. If $|X|+|N_{H_{2}}(s)-N_{H_{2}}(v_{t})|>|N_{F}(s)\cap N_{H_{2}}(v_{t})|,$ then by Lemma~\ref{lem:existsAnExch} and Lemma~\ref{lem:reduce}, for some $x\in X$, there exists a simplified $v_{t}x$ and $xs$ exchange that when exchanged creates another $(G',F',t')\in \mathcal{G}$. Note that since the exchange was simplified, $M-\{sx\}$ is a matching of $G'-E(F')$. If $x\notin V(M)$, then $M+\{v_{t}x\}$ is a matching of $G'-E(F')$ that contradicts \ref{highDegreeMatching_choice:1}. If $x\in V(M)$, then since $\overline{x}_{M}\in Q$ by Claim~\ref{cl:inA}, we see that $M-\{x\overline{x}_{M}\}+\{v_{t}x\}$ is a matching of $G'-E(F')$ that violates \ref{highDegreeMatching_choice:2}. Thus, we are left with the case \begin{equation}\label{eq:mainproof}|N_{F}(s)\cap N_{H_{2}}(v_{t})|\geq |X|+|N_{H_{2}}(s)-N_{H_{2}}(v_{t})|\geq |X|+|N_{H_{2}}(s)\cap N_{F}(v_{t})|.
\end{equation}
Considering (\ref{eq:mainproof}) we may deduce that 
\begin{align}|N_{F}(s)-N_{H_{1}}(v_{t})|&=|N_{F}(s)\cap N_{F}(v_{t})|+|N_{F}(s)\cap N_{H_{2}}(v_{t})|\notag \\
&\geq |N_{F}(s)\cap N_{F}(v_{t})|+|X|+|N_{H_{2}}(s)\cap N_{F}(v_{t})|\notag\\
&\geq |N_{F}(s)|-|N_{H_{1}}(s)\cap N_{F}(v_{t})|+|X|\notag \\
&=r-|N_{H_{1}}(s)\cap N_{F}(v_{t})|+|X|.\label{eq:mainproof2}
\end{align}

Observe that \begin{equation}|N_{H_{1}}(s)- N_{H_{1}}(v_{t})|=|N_{H_{1}}(s)\cap N_{F}(v_{t})|+|N_{H_{1}}(s)\cap N_{H_{2}}(v_{t})|, \label{eq:mainproof3}
\end{equation}
and \begin{equation}|N_{H_{1}}(s)\cap N_{H_{2}}(v_{t})|+|X|\geq |(N_{H_{1}}(s)\cap D')-N_{H_{1}}(v_{t})|\geq d_{n}-r+1. \label{eq:mainproof4}
\end{equation}
From (\ref{eq:mainproof3}) and (\ref{eq:mainproof4}), we deduce that \begin{equation}|N_{H_{1}}(s)- N_{H_{1}}(v_{t})|\geq |N_{H_{1}}(s)\cap N_{F}(v_{t})|+ d_{n}-r+1-|X|. \label{eq:mainproof5}
\end{equation}

Thus, combining (\ref{eq:mainproof2}) and (\ref{eq:mainproof5}) we see that \begin{equation}|N_{H_{1}}(s)- N_{H_{1}}(v_{t})|+|N_{F}(s)-N_{H_{1}}(v_{t})|\geq d_{n}+1. \label{eq:mainproof6}
\end{equation}

Since $deg_{H_{1}}(v_{t})\geq d_{1}-d_{n}$ and Claim~\ref{cl:inA} states that $N_{H_{1}}(v_{t})\subseteq N_{G}(s)$, we can use (\ref{eq:mainproof6}) to show our final contradiction
\begin{align}
d_{1}\geq deg_{G}(s)&\geq deg_{H_{1}}(v_{t})+|N_{H_{1}}(s)- N_{H_{1}}(v_{t})|+|N_{F}(s)-N_{H_{1}}(v_{t})|\notag\\
&\geq d_{1}-d_{n}+d_{n}+1= d_{1}+1. \qedhere
\end{align}
\end{proof}

\section{Proof of Theorem~\lowercase{\ref{thm:middegreebound}}}\label{sec:middegreebound}

\begin{proof}We assume every $G\in \mathcal{R}(\pi)$ has vertex set $V=\{v_{1},\ldots,v_{n}\}$ such that $deg_{G}(v_{i})=d_{i}$. Suppose $\big\lceil\frac{n+5-k}{2}\big\rceil>d_{\frac{n}{2}+1}$.  Since $\mathcal{D}_{k}(\pi)$ is graphic, Kundu's $k$-factor Theorem says $\pi$ has a realization with a $k$-factor. Moreover, by Theorem~\ref{thm:main-split} and Petersen's $2$-factor theorem, $\pi$ has a realization with a $k$-factor whose edges can be partitioned into $1$-factors and $2$-factors. We let $r$ be the largest natural number such that $k+r$ is even and there is a $G\in \mathcal{R}(\pi)$ with a $k$-factor $F$ whose edges can be partitioned into graphs $F_{1},\dots,F_{\frac{k+r}{2}}$ such that $F_{i}$ is a $1$-factor when $i\leq r$ and a $2$-factor when $i>r$. By contradiction we assume $r\leq k-2$. Note that for $i>r$, every $F_{i}$ must have at least two odd cycles as components. Otherwise, we could split $F_{i}$ into two $1$-factors contradicting our choice of $r$. 

We let $(H,H_{1},H_{2},H_{3},\ldots, H_{q})$ correspond to an edge coloring of $K_{n}$ with the natural numbers $\{1,\ldots,q\}$ such that $H=K_{n}-E(H_{2})$ and $H_{i}$ is the subgraph induced by the edges colored $i$. We let $\mathcal{G}$ be the set of such tuples where $H_{1}\in \mathcal{R}(G-E(F))$, $H_{2}\in \mathcal{R}(\overline{G})$, and $H_{i}\in \mathcal{R}(F_{i-2})$ for all $3\leq i\leq q$. We further consider the subset $\mathcal{G}'\subseteq \mathcal{G}$ to be all tuples such that the number of cycles in $H_{q}$ is minimized.

For the rest of this proof, we will assume all color exchanges are simplified.

\begin{claim}\label{cl:contstructLargeCycle}For $(H,H_{1},\ldots, H_{q})\in \mathcal{G}'$, let $A=a_{0}\ldots a_{|A|-1}a_{0}$ and $B=b_{0}\ldots b_{|B|-1}b_{0}$ be distinct cycles of $H_{q}$. For some $a_{i}$ and $b_{j}$, if $L=\{a_{i}x_{0},x_{0}b_{j},\ldots, a_{i}x_{l-1},x_{l-1}b_{j}\}$ is a color exchange, then $\mathcal{X}(L)\cap (V(A)\cup V(B))=\emptyset$.
\begin{proof}Suppose the claim is false and there is an edge $a_{i}x_{s}$ of $A$ and an edge $x_{s+1}b_{j}$ of $B$. After, performing the exchange we denote the resulting tuple as $(H',H'_{1},\ldots, H'_{q})$. Since $L$ is simplified $a_{i}x_{s}$ and $x_{s+1}b_{j}$ are the only edges of $H_{q}$ exchanged. Without loss of generality we assume $x_{s}=a_{i'}$ and $x_{s+1}=b_{j'}$ with $i'= i+1 \mod |A|$ and $j'= j+1 \mod |B|$. With this we can see
$a_{i'}\ldots a_{|A|-1}a_{0}\ldots a_{i}b_{j'}\dots b_{|B|-1}b_{0}\ldots b_{j}a_{i'}$ is a cycle in $H'_{q}$ that combines the vertices of $A$ and $B$ and leaves all other cycles alone. However, this implies the contradiction that $H'_{q}$ has fewer cycles than $H_{q}$. 
\end{proof}
\end{claim}

We choose an arbitrary $(H,H_{1},\ldots, H_{q})\in \mathcal{G}'$. We let $f$ be the largest index such that $v_{f}v^{-}_{f}$ is an edge of $H_{q}$ with $deg_{H}(v_{f})\geq deg_{H}(v^{-}_{f})$ and denote the cycle containing $v_{f}v^{-}_{f}$ by $A$. Since $H_{q}$ does not have a $1$-factor, there is an odd cycle $C$ that is distinct from $A$. Therefore, there is a smallest $t$ such that $v_{t}\in V(C)$ and the two neighbors $v^{-}_{t}$ and $v^{+}_{t}$ of $v_{t}$ along $C$ satisfy $deg_{G}(v^{-}_{t})\leq deg_{G}(v_{t})\leq deg_{G}(v^{+}_{t})$. By our choice of $v_{f}$, we know that $t<f$. There is an odd cycle $D$ in $H_{q}$ that is not $C$. Like $v_{t}$ in $C$, $D$ has vertices $v_{s}$, $v^{+}_{s}$, and $v^{-}_{s}$ such that $deg_{H_{q}}(v^{-}_{s})\leq deg_{H_{q}}(v_{s}) \leq deg_{H_{q}}(v^{+}_{s})$. By our choice of $v_{f}$, we know that $f\geq s$.

Since we maximized our selection of $f$, we know that $\{v_{f+1},\ldots,v_{n}\}$ is an independent set in $H_{q}$. Thus, $e_{H_{q}}(\{v_{1},\ldots,v_{f}\},\{v_{f+1},\ldots,v_{n}\})= 2(n-f)$.  Note that $\{v_{s},v^{+}_{s},v_{t},v^{+}_{t}\}\subseteq \{v_{1},\ldots,v_{f}\}$ and $v_{s}v^{+}_{s}$ and $v_{t}v^{+}_{t}$ are edges of $H_{q}$. Therefore, \[e_{H_{q}}(\{v_{1},\ldots,v_{f}\},\{v_{f+1},\ldots,v_{n}\})\leq 2f-|\{v_{s},v^{+}_{s},v_{t},v^{+}_{t}\}|=2(f-2).\] Combining the two bounds of $e_{H_{q}}(\{v_{1},\ldots,v_{f}\},\{v_{f+1},\ldots,v_{n}\})$ and solving for $f$, we deduce that $f\geq \frac{n}{2}+1$.

\begin{claim}\label{cl:colorfulEdge}\[\{v_{f}v^{+}_{t},v_{f}v_{t},v^{-}_{f}v^{+}_{t},v^{-}_{f}v_{t}\}\subseteq \bigcup_{3\leq i<q} E(H_{i}).\]
\begin{proof}Note that $deg_{H}(v^{-}_{f})\leq deg_{H}(v_{f})\leq deg_{H}(v_{t})\leq deg_{H}(v^{+}_{t})$. If $v_{f}v_{t}$ is an edge of $H_{1}$, then we deduce a contradiction to Claim~\ref{cl:contstructLargeCycle} since $v^{+}_{t}v_{t}$ and $v_{t}v_{f}$ can be exchanged. If $v_{f}v_{t}$ is an edge of $H_{2}$, then we deduce a contradiction to Claim~\ref{cl:contstructLargeCycle} since $v^{-}_{f}v_{f}$ and $v_{f}v_{t}$ can be exchanged. Similar arguments prove the claim for $v_{f}v^{+}_{t}$, $v^{-}_{f}v^{+}_{t}$, and $v^{-}_{f}v_{t}$.
\end{proof}
\end{claim}

\begin{claim}\label{cl:H1neighbors}$N_{H_{1}}(v_{t})\subseteq N_{H}(v_{f})\cap \{v_{1},\ldots,v_{f}\}$. 
\begin{proof}Let $v_{i}\in N_{H_{1}}(v_{t})$.
Suppose $v_{i}\in V(A)$. If $i>t$, then by Lemma~\ref{lem:largedegree} we deduce a contradiction to Claim~\ref{cl:contstructLargeCycle} since $v^{+}_{t}v_{t}$ and $v_{t}v_{i}$ can be exchanged. Thus, $f>t>i$. If $v_{i}\in N_{H_{2}}(v_{f})$, then by Claim~\ref{cl:contstructLargeCycle} and Claim~\ref{cl:colorfulEdge}, we can exchange the edges $v_{t}v_{f}$ and $v_{f}v_{i}$ to create a tuple $(H'_{1},\ldots,H'_{q})\in \mathcal{G}'$ such that $H_{q}=H'_{q}$.  However, we see a contradiction to Claim~\ref{cl:contstructLargeCycle} since $v_{f}^{-}v_{f}$ and $v_{f}v_{t}$ can be exchanged with respect to this new tuple. Thus, $v_{i}\subseteq N_{H}(v_{f})\cap \{v_{1},\ldots,v_{t-1}\}\subset N_{H}(v_{f})\cap \{v_{1},\ldots,v_{f}\}$ when $v_{i}\in V(A)$.

Suppose $v_{i}\notin V(A)$. If $i>f$, then by Claim~\ref{cl:contstructLargeCycle} and Claim~\ref{cl:colorfulEdge}, we can exchange the edges $v_{f}v_{t}$ and $v_{t}v_{i}$ to create a tuple $(H'_{1},\ldots,H'_{q})\in \mathcal{G}'$ such that $H_{q}=H'_{q}$.  However, we have a contradiction to Claim~\ref{cl:contstructLargeCycle} since $v_{t}^{+}v_{t}$ and $v_{t}v_{f}$ can be exchanged with respect to this new tuple. Thus, we are left with the case $i<f$ and $v_{i}\in N_{H_{2}}(v_{f})$. Here we have a contradiction to Claim~\ref{cl:contstructLargeCycle} since $v_{f}^{-}v_{f}$ and $v_{f}v_{i}$ can be exchanged. Thus, $v_{i}\subseteq N_{H}(v_{f})\cap \{v_{1},\ldots,v_{f}\}$ when $v_{i}\notin V(A)$.
\end{proof}
\end{claim}

Let $\overline{\mathcal{D}_{k}(\pi)}=(\overline{q}_{1},\dots, \overline{q}_{n})$ where $\overline{q}_{i}=n-1-d_{n+1-i}+k$. 

Since $\{v^{+}_{t},v_{t},v_{f}\}\subseteq \{v_{1},\ldots, v_{f}\}$ and not in $N_{H_{1}}(v_{t})$, Claim~\ref{cl:H1neighbors} implies that $|N_{H_{1}}(v_{t})|\leq f-3$. For $v_{i}\in \{v_{f-2},\ldots,v_{n}\}$, this implies \[\overline{q}_{n+1-i}=n-1-d_{i}+k\geq n-1-|N_{H_{1}}(v_{t})|\geq n+2-f.\] Thus, $m(\overline{D_{k}(\pi)})\geq n+3-f$, and therefore, $f\geq \max\{\frac{n}{2}+1, n+3-m(\overline{D_{k}(\pi)})\}$.

We now turn to finding a lower bound for $d_{f}$.

\begin{claim}\label{cl:colorNeighbors}If $v_{i}v_{j}\in E(H_{q})$ such that $v_{i}\in N_{H_{2}}(v_{t})$, then $v_{j}\in N_{H}(v_{f})\cap \{v_{1},\ldots,v_{f}\}$.
\begin{proof}Suppose $v_{j}\notin N_{H}(v_{f})\cap \{v_{1},\ldots,v_{f}\}$. We first assume $v_{i}\notin V(C)$. By Claim~\ref{cl:contstructLargeCycle} $v_{j}v_{i}$ and $v_{i}v_{t}$ cannot be exchanged. Therefore, by Lemma~\ref{lem:largedegree} $j<t$. If $v_{j}\in N_{H_{2}}(v_{f})$, then we can exchange the edges $v_{t}v_{f}$ and $v_{f}v_{j}$ to create a tuple $(H'_{1},\ldots,H'_{q})\in \mathcal{G}'$ such that $H_{q}=H'_{q}$.  However, we have a contradiction to Claim~\ref{cl:contstructLargeCycle} since $v_{f}^{-}v_{f}$ and $v_{f}v_{t}$ can be exchanged with respect to this new tuple. Thus, $v_{j}\in N_{H}(v_{f})\cap \{v_{1},\ldots,v_{f}\}$ when $v_{i}\notin V(C)$.

We are left with the case $v_{i}\in V(C)$. If $i<f$, then we can exchange the edges $v_{f}v_{t}$ and $v_{t}v_{i}$ to create a tuple $(H'_{1},\ldots,H'_{q})\in \mathcal{G}'$ such that $H_{q}=H'_{q}$.  However, we see a contradiction to Claim~\ref{cl:contstructLargeCycle} since $v_{f}^{-}v_{f}$ and $v_{f}v_{t}$ can be exchanged with respect to this new tuple. Thus, $i>f$, and therefore, $j<f$ since we chose the largest such $f$. By Claim~\ref{cl:contstructLargeCycle} we know that $v_{i}v_{j}$ and $v_{j}v_{f}$ cannot be exchanged. Thus, by Lemma~\ref{lem:largedegree} $v_{j}\in N_{H}(v_{f})$, and therefore, $v_{j}\in N_{H}(v_{f})\cap \{v_{1},\ldots,v_{f}\}$.
\end{proof}
\end{claim}

We let $W=N_{H}(v_{f})\cap \{v_{1},\ldots,v_{f}\}$. Claim~\ref{cl:H1neighbors} states that $N_{H_{1}}(v_{t})\subseteq W$, and we know that $v_{t}$ is adjacent in $H-E(H_{1})$ to $v_{f}$, $v_{f}^{-}$, and $v_{t}^{+}$. Thus, \[d_{f}=|N_{H}(v_{f})|\geq |N_{H_{1}}(v_{t})\cup \{v_{f}^{-},v^{+}_{t},v_{t}\}|\geq d_{t}-k+3.\] Since Claim~\ref{cl:colorNeighbors} states $N_{H_{q}}(v_{i})\subseteq W$ for all $v_{i}\in N_{H_{2}}(v_{t})$, we deduce that \[e_{H_{q}}(N_{H_{2}}(v_{t}),W)\geq 2|N_{H_{2}}(v_{t})|=2(n-1-d_{t}).\] On the other hand, $v^{+}_{t}$ and $v^{+}_{s}$ are in $W$ and each of them are adjacent in $H_{q}$ to at least one vertex not in $N_{H_{2}}(v_{t})$, and $v_{t}$ is adjacent in $H_{q}$ to two vertices not in $N_{H_{2}}(v_{t})$. Note that it is possible for $v_{s}=v_{f}$. Thus, we left $v_{s}v^{-}_{s}$ out of our calculation. Therefore, \[e_{H_{q}}(N_{H_{2}}(v_{t}), W)\leq 2|W|-4.\] Combining we see that \[2|W|-4\geq e_{H_{q}}(N_{H_{2}}(v_{t}), W)\geq 2(n-1-d_{t}).\] Solving for $|W|$ we deduce that \[|W|\geq n+1-d_{t}.\] Since $v^{-}_{f}\notin W$, we see that $|W|\leq d_{f}-1$. Therefore, \[\max\{ n+2-d_{t},d_{t}-k+3\}\leq d_{f}.\] Letting $n+2-d_{t}=d_{t}-k+3$ we deduce that $\max\{ n+2-d_{t},d_{t}-k+3\}$ is minimized when $2d_{t}-(n-1)=k$. Thus, \[\bigg\lceil\frac{n+5-k}{2}\bigg\rceil\leq d_{f}.\] Using the lower bound on $f$, we see the contradiction \[\bigg\lceil\frac{n+5-k}{2}\bigg\rceil\leq d_{f}\leq d_{\max\{\frac{n}{2}+1, n+3-m(\overline{D_{k}(\pi)})\}}.\]

We let $l = \max\{\frac{n}{2}+1, n+3-m(\pi)\}$, and suppose \[\bigg\lceil\frac{n+5-k}{2}\bigg\rceil\leq \overline{q}_{l}.\] Therefore, $\overline{q}_{l}=n-1-d_{n+1-l}+k\geq \bigg\lceil\frac{n+5-k}{2}\bigg\rceil$. Solving for $d_{n+1- l}$, and realizing $n+1-l=\min\{\frac{n}{2},m(\pi)-2\}$ we deduce the contradiction \[d_{\min\{\frac{n}{2},m(\pi)-2\}}\leq n+k-1-\bigg\lceil\frac{n+5-k}{2}\bigg\rceil=\bigg\lceil\frac{n+3k-8}{2}\bigg\rceil.\] Thus, \[ \overline{q}_{l}< \bigg\lceil\frac{n+5-k}{2}\bigg\rceil,\] and by the first part of this theorem, some realization of $\overline{D_{k}(\pi)}$ has $k$-edge-disjoint $1$-factors. Therefore, some realization of $\pi$ has $k$ edge-disjoint $1$-factors.
\end{proof}

\section{Proof of theorem~\lowercase{\ref{thm:main-split}}}\label{sec:main-split}
\begin{proof}Let $\mathcal{G}$ be the set of tuples of the form $(G,F,r,F_{0})$ where $G\in \mathcal{R}(\pi)$, $F$ is a $k$-factor of $G$, $F_{0}$ is a spanning $(k-r)$-factor of $F$, and $E(F-E(F_{0}))$ can be partitioned into $r$ edge-disjoint $1$-factors.
\begin{enumerate}[label=(C\arabic*),ref=(C\arabic*)]
    \item\label{highKMatching_choice:1} We choose a $(G,F,r, F_{0})\in \mathcal{G}$ such that $r$ is maximized, and
    \item\label{highKMatching_choice:2} subject to $\ref{highKMatching_choice:1}$, we minimize $\text{def}(F_{0})$. 
\end{enumerate}

Let $k'=k-r$. It is enough to show that $r\geq \big \lceil\frac{k'+1}{2}\big\rceil+2$ since this would imply $r\geq \min\{\lceil\frac{k+5}{3}\big\rceil,k\}$. Let $M_{1},\ldots, M_{r}$ be a partition of $E(F-E(F_{0}))$ into $r$ edge-disjoint $1$-factors. 

Let $H_{1}=G-E(F)$, $H_{2}=\overline{G}$, and $H_{3}=F_{0}$. Note that $H_{1},H_{2}, H_{3},\ldots, H_{r+2}$ where $H_{i}=M_{i-3}$ for $i\geq 4$ represent a coloring of the edges of $K_{n}$. Thus, any edge exchange involving any $H_{i}$ in $K_{n}$ corresponds to some $(G',F',r',F'_{0})\in \mathcal{G}$.

Let $A$, $C$, and $D$ be a Gallai-Edmonds Decomposition of $F_{0}$. Since $F_{0}$ does not have a $1$-factor we know that $D$ is not empty. Let $\mathcal{D}=\{D_{1},\ldots,D_{|\mathcal{D}|}\}$ be the components of $F_{0}[D]$.  

\begin{claim}\label{cl:stars}Suppose $u\in V(D_{j})$ and $x\in V(D_{i})$ such that there exists a maximum matching $M$ of $F_{0}$ that misses both $u$ and $x$. If $uy\in E(H_{f})$ for some $y\in N_{D_{i}}(x)$, then $xv\notin E(H_{f})$ for all $v\in N_{D_{j}}(u)$.
\begin{proof}If there does exists such a $v\in N_{D_{j}}(u)$, then we may exchange the edges $uv$ and $xy$ of $F_{0}$ with the edges $xv$ and $uy$ of $H_{f}$ to find a $(G',F',r,F'_{0})\in \mathcal{G}$ such that $M+\{ux\}$ is a matching in $F'_{0}$. This contradicts \ref{highDegreeMatching_choice:2}, since $\text{def}(F'_{0})\leq \text{def}(F_{0})-2$.
\end{proof}
\end{claim}

\begin{claim}\label{cl:starsObservation}If $v_{2}\in V(D_{j})$ and $u_{1}\in V(D_{i})$ such that there exists a maximum matching $M$ of $F_{0}$ that misses both $v_{2}$ and $u_{1}$, then \[r\geq |N_{F-E(F_{0})}(u_{1})\cap N_{D_{j}}[v_{2}]|+|N_{F-E(F_{0})}(v_{2})\cap N_{D_{i}}[u_{1}]|-1.\]
\begin{proof}
    For any $v_{1}\in N_{F-E(F_{0})}(u_{1})\cap N_{D_{j}}[v_{2}]$ and $u_{2}\in N_{F-E(F_{0})}(v_{2})\cap N_{D_{i}}[u_{1}]$, Claim~\ref{cl:stars} implies $\{v_{1}u_{2},v_{1}u_{1},v_{2}u_{2}\}$ are in different $1$-factors. Thus, simply adding up all of the distinct edges of $F-E(F_{0})$ between $N_{F-E(F_{0})}(u_{1})\cap N_{D_{j}}[v_{2}]$ and $N_{F-E(F_{0})}(v_{2})\cap N_{D_{i}}[u_{1}]$ establishes the claim.
\end{proof}
\end{claim}

We let $L_{i}$ denote the set of all vertices in $D_{i}$ such that for every $u\in L_{i}$, there are edges $uu^{-}$ and $uu^{+}$ of $D_{i}$ with $deg_{G}(u^{+})\geq deg_{G}(u)\geq deg_{G}(u^{-})$.

Lov{\'a}sz \cite{Lovasz1972} showed that every non-trivial factor critical graph has an odd cycle. Since Theorem~\ref{thm:GallieEdmonds}~\ref{GEST:2} states that every component in $\mathcal{D}$ is factor critical, each non-trivial component must have an odd cycle. Moreover, for any odd cycle $C_{i}$ in $D_{i}$, there must be distinct edges $\{uu^{+},uu^{-}\}$ of $C_{i}$ such that $deg_{G}(u^{+})\geq deg_{G}(u)\geq deg_{G}(u^{-})$. Thus, $L_{i}$ is not empty when $D_{i}$ is not trivial.

By \ref{GEST:1} of Theorem~\ref{thm:GallieEdmonds}, $|\mathcal{D}|\geq |A|+2$. Since $k'(|A|+2)\leq k'|D|\leq e_{F_{0}}(D,A)\leq k'|A|$, there must be at least two non-trivial components in $\mathcal{D}$. Thus, $\bigcup_{i=1}^{|\mathcal{D}|} L_{i}$ has at least two vertices.

Let $z_{i}\in L_{i}$, and let $M$ be a maximum matching of $F_{0}$ that misses $z_{i}$. We let $D'(z_{i},M)\subseteq D-V(D_{i})$ be the largest set such that for every $u\in D'(z_{i},M)$ there is a maximum matching $M_{u}$ of $F_{0}$ that misses both $u$ and some $v\in D_{i}$ such that $E(M_{u}-D'(z_{i},M))=E(M-D'(z_{i},M))$. Since every component in $\mathcal{D}$ is factor critical and $M$ restricts to a near matching on each of them, we may conclude that for $D_{j}\in \mathcal{D}$, if $V(D_{j})\cap D'(z_{i},M)\neq \emptyset$, then $V(D_{j})\subseteq D'(z_{i},M)$. We let $\mathcal{D}'(z_{i},M)$ be the components of $F_{0}[D'(z_{i},M)]$.

\begin{claim}\label{cl:NumberofComp} $|\mathcal{D}'(z_{i},M)|\geq |N_{F_{0}}(D'(z_{i},M))\cap A|+1$ for a given $z_{i}\in L_{i}$ and a given maximum matching $M$ of $F_{0}$ that misses $z_{i}$.
\begin{proof}
     Let $S=N_{F_{0}}(D'(z_{i},M))\cap A$. Suppose $\overline{w}_{M}\notin D'(z_{i},M)$ for some $w\in S$. By the definition of $D'(z_{i},M)$ for every $v\in D'(z_{i},M)$ there is a $M_{v}$ such that $w\overline{w}_{M}\in E(M_{v})$, and therefore, $E(M_{v}-(D'(z_{i},M)\cup \{\overline{w}_{M}\}))=E(M-(D'(z_{i},M)\cup \{\overline{w}_{M}\}))$. Let $u\in N_{F_{0}}(w)\cap D'(z_{i},M)$. By the definition of $D'(z_{i},M)$, there is a maximum matching $M_{u}$ of $F_{0}$ that misses $u$ and some $v\in V(D_{i})$. Since $w\overline{w}_{M}\in E(M_{u})$, the matching $M_{u}-\{w\overline{w}_{M}\}+\{wu\}$ misses $\overline{w}_{M}$ and $z_{i}$. Moreover, $E((M_{u}-\{w\overline{w}_{M}\}+\{wu\})-(D'(z_{i},M)\cup \{\overline{w}_{M}\}))=E(M-(D'(z_{i},M)\cup \{\overline{w}_{M}\}))$. Thus, we have shown that for every $u\in D'(z_{i},M)\cup \{\overline{w}_{M}\}$ there is some $v\in V(D_{i})$ and a maximum matching $M_{u}$ that misses $u$ and $v$ such that $E(M_{u}-(D'(z_{i},M)\cup \{\overline{w}_{M}\}))=E(M-(D'(z_{i},M)\cup \{\overline{w}_{M}\}))$. However, this is a contradiction since $|D'(z_{i},M)\cup \{\overline{w}_{M}\}|>|D'(z_{i},M)|$. Thus, $\overline{w}\in D'(z_{i},M)$ for all $w\in S$, and therefore, \ref{GEST:1} of Theorem~\ref{thm:GallieEdmonds} implies $|\mathcal{D}'(z_{i},M)|\geq |S|+1$.
\end{proof}
\end{claim}

For $u\in L_{i}$ and $v\in V$, we let $W^{+}_{i}(u,v)$ be all $w\in N_{D_{i}}[u]$ such that $deg_{G}(w)\geq deg_{G}(v)$. Similarly, we let $W^{-}_{i}(u,v)$ be all $w\in N_{D_{i}}[u]$ such that $deg_{G}(w)\leq deg_{G}(v)$.

\begin{claim}\label{cl:forbidden} Suppose $v_{2}\in L_{j}$ and $u_{1}\in L_{i}$ such that there exists a maximum matching $M$ of $F_{0}$ that misses both $v_{2}$ and $u_{1}$. If $v_{1}\in W_{j}^{+}(v_{2},u_{1})-\{v_{2}\}$ and $u_{2}\in W_{i}^{-}(u_{1},v_{2})-\{u_{1}\}$, then $v_{1}u_{2}$ is not an edge of $H_{1}$ nor $H_{2}$. Moreover, if $deg_{G}(v_{2})\geq deg_{G}(u_{1})$, then the edges $\{u_{1}v_{1}, u_{1}v_{2}, u_{2}v_{1}, u_{2}v_{2}\}$ are in distinct $H_{s}$ for $s\geq 3$.
\begin{proof}Since $D_{j}$ and $D_{i}$ are factor critical and $M$ restricts to a near perfect matching on $D_{i}$ and $D_{j}$, we can find a maximum matching $M'$ of $F_{0}$ that misses both $v_{1}$ and $u_{2}$.

Suppose $v_{1}u_{2}\in E(H_{1})$. Since $deg_{G}(v_{2})\geq deg_{G}(u_{2})$, Lemmas \ref{lem:reduce} and \ref{lem:largedegree} imply there exists a simplified $v_{2}v_{1}$ and $v_{1}u_{2}$ edge-exchange that when exchanged creates another $(G',F',r,F'_{0})\in \mathcal{G}$ with a larger matching $M'+\{v_{1}u_{2}\}$ that contradicts \ref{highDegreeMatching_choice:2}. Thus, $v_{1}u_{2}\notin E(H_{1})$.

Suppose $v_{1}u_{2}\in E(H_{2})$. Since $deg_{G}(v_{1})\geq deg_{G}(u_{1})$, Lemmas \ref{lem:reduce} and \ref{lem:largedegree} imply there exists a simplified $u_{2}u_{1}$ and $v_{1}u_{2}$ edge-exchange that when exchanged creates another $(G',F',r,F'_{0})\in \mathcal{G}$ with a larger matching $M'+\{v_{1}u_{2}\}$ that contradicts \ref{highDegreeMatching_choice:2}. Thus, $v_{1}u_{2}\notin E(H_{2})$.

We may now apply the first part of the claim to the second. Suppose $deg_{G}(v_{2})\geq deg_{G}(u_{1})$. Since $v_{2}\in W_{j}^{+}(v_{1},u_{1})-\{v_{1}\}$ and $u_{2}\in W_{i}^{-}(u_{1},v_{1})-\{u_{1}\}$, $v_{2}u_{2}$ is not an edge of $H_{1}$ nor $H_{2}$. Since $v_{1}\in W_{j}^{+}(v_{2},u_{2})-\{v_{2}\}$ and $u_{1}\in W_{i}^{-}(u_{2},v_{2})-\{u_{2}\}$, $v_{1}u_{1}\notin E(H_{1})\cup E(H_{2})$. Since $v_{2}\in W_{j}^{+}(v_{1},u_{1})-\{v_{1}\}$ and $u_{1}\in W_{i}^{-}(u_{2},v_{1})-\{u_{2}\}$, $v_{2}u_{1}\notin E(H_{1})\cup E(H_{2})$.  Moreover, Claim~\ref{cl:stars} implies the edges $\{u_{1}v_{1}, u_{1}v_{2}, u_{2}v_{1}, u_{2}v_{2}\}$ are in distinct $H_{s}$ for $s\geq 3$.
\end{proof}
\end{claim}

\begin{claim}\label{cl:forbidden2}Suppose $v_{2}\in L_{j}$ and $u_{1}\in L_{i}$  such that there exists a maximum matching $M$ of $F_{0}$ that misses both $v_{2}$ and $u_{1}$. If $v_{1}\in W_{j}^{+}(v_{2},u_{1})$ and $u_{2}\in W_{i}^{-}(u_{1},v_{2})$, then $v_{1}u_{2}$ is not an edge of $H_{1}$ nor $H_{2}$.
\begin{proof}If $v_{2}\notin W_{j}^{+}(v_{2},u_{1})$, then $u_{1}\notin W_{i}^{-}(u_{1},v_{2})$, and therefore, $v_{1}\in W_{j}^{+}(v_{2},u_{1})-\{v_{2}\}$ and $u_{2}\in W_{j}^{-}(u_{1},v_{2})-\{u_{1}\}$. Thus, the Claim follows from Claim~\ref{cl:forbidden}. If $v_{2}\in W_{j}^{+}(v_{2},u_{1})$, then $u_{1}\in W_{i}^{-}(v_{2},u_{1})$. In this case, since $v_{2}\in L_{j}$ and $u_{1}\in L_{i}$, $|W^{+}(v_{2},u_{1})|\geq 2$ and $|W^{-}(v_{2},u_{1})|\geq 2$. Thus, there are vertices $v_{3}$ and $u_{3}$ such that $v_{1}\in W_{j}^{+}(v_{3},u_{1})-\{v_{3}\}$ and $u_{2}\in W_{i}^{-}(u_{3},v_{2})-\{u_{3}\}$. Thus, the Claim follows from Claim~\ref{cl:forbidden}.
\end{proof}
\end{claim}

We choose a $z_{1}\in \bigcup_{i=1}^{|\mathcal{D}|} L_{i}$ such that $e_{F_{0}}(z_{1},A)$ is minimized, and subject to that, we give preference to one that maximizes $deg_{G}(z_{1})$. Without loss of generality we may assume $z_{1}\in L_{1}$, and therefore, $D_{1}$ is not trivial.

Let $M$ be a maximum matching of $F_{0}$ that misses $z_{1}$. We let $D'=D'(z_{1},M)$, Let $S=N_{F_{0}}(D')\cap A$, and $\mathcal{D}'=\mathcal{D}'(z_{1},M)$. Without loss of generality we assume $\mathcal{D}'=\{D_{2},\ldots,D_{|\mathcal{D}'|+1}$\}.

Since $|S|k'\geq e_{F_{0}}(D',S)$ and there are at least $|S|+1$ components in $\mathcal{D}'$, we may conclude that there is at least one component in $\mathcal{D}'$ that is adjacent in $F_{0}$ to at most $k'-1$ vertices in $S$. Without loss of generality we may assume $e_{F_{0}}(D_{i},S)\leq k-1$ for $D_{i}\in \{D_{2},\ldots, D_{t}\}$ where $2\leq t\leq |\mathcal{D}'|+1$. Clearly, the components in $\{D_{2},\ldots, D_{t}\}$ are non-trivial. 

Thus, \[k'|S|\geq \sum_{i=2}^{|\mathcal{D}'|}e_{F_{0}}(D_{i},S)\geq k'(|S|+1-(t-1))+\sum_{i=2}^{t}e_{F_{0}}(D_{i},S).\] Therefore, \begin{equation}\label{eq:forbidden}\sum_{i=2}^{t}e_{F_{0}}(D_{i},S)\leq k'(t-2).
\end{equation}

Our proof will contain several common arguments, so we will address them first in the next two claims.
\begin{claim}\label{cl:join}For $y\in L_{1}$ and $u\in L_{i}$ where $D_{i}\in \{D_{2},\ldots, D_{t}\}$, every vertex in $W^{+}_{i}(u,z_{1})\cap W^{+}_{i}(u,y)$ is adjacent in $F-E(F_{0})$ to every vertex in $W^{-}_{1}(y,u)\cup W^{-}_{1}(z_{1},u)$, and every vertex in $W^{-}_{i}(u,z_{1})\cap W^{-}_{i}(u,y)$ is adjacent in $F-E(F_{0})$ to every vertex in $W^{+}_{1}(y,u)\cup W^{+}_{1}(z_{1},u)$.
\begin{proof}This follows by repeatedly applying Claim~\ref{cl:forbidden2} with respect to $M$. 
\end{proof}
\end{claim}

\begin{claim}\label{cl:joinDistinct}For $y\in L_{1}$ and $\{u,v\}\in L_{i}$ where $D_{i}\in \{D_{2},\ldots, D_{t}\}$ and $u$ could equal $v$, the following holds:
\begin{enumerate}
    \item \label{it1} If $W^{+}_{i}(u,z_{1})\cap W^{+}_{i}(u,y)\neq \emptyset$, then \[r\geq \max\{|W^{-}_{1}(z_{1},u)\cup W^{-}_{1}(y,u)|, |W^{-}_{1}(z_{1},u)\cap W^{-}_{1}(y,u)|+|\{z_{1},y\}|\}.\]
    \item \label{it2} If $W^{-}_{i}(u,z_{1})\cap W^{-}_{i}(u,y)\neq \emptyset$, then \[r\geq \max\{|W^{+}_{1}(z_{1},u)\cup W^{+}_{1}(y,u)|, |W^{+}_{1}(z_{1},u)\cap W^{+}_{1}(y,u)|+|\{z_{1},y\}|\}.\]
    \item \label{it3} If $W^{+}_{1}(y,u)\cap W^{+}_{1}(y,v)\neq \emptyset$, then \[r\geq \max\{|W^{-}_{i}(u,y)\cup W^{-}_{i}(v,y)|, |W^{-}_{i}(u,y)\cap W^{-}_{i}(v,y)|+|\{u,v\}|\}.\]
    \item \label{it4} If $W^{-}_{1}(y,u)\cap W^{-}_{1}(y,v)\neq \emptyset$, then \[r\geq \max\{|W^{+}_{i}(u,y)\cup W^{+}_{i}(v,y)|, |W^{+}_{i}(u,y)\cap W^{+}_{i}(v,y)|+|\{u,v\}|\}.\]
    \item \label{it9} If $deg_{G}(z_{1})\leq deg_{G}(u)$, then $r\geq |W^{+}_{i}(u,z_{1})|+|W^{-}_{1}(z_{1},u)|-1$
    \item \label{it10} If $deg_{G}(z_{1})\geq deg_{G}(u)$, then $r\geq |W^{-}_{i}(u,z_{1})|+|W^{+}_{1}(z_{1},u)|-1$.
\end{enumerate}
\begin{proof}
    Suppose there is a $w\in W^{+}_{i}(u,z_{1})\cap W^{+}_{i}(u,y)$. By Claim~\ref{cl:join} $w$ is adjacent in $F-E(F_{0})$ to every vertex in $W^{-}_{1}(z_{1},u)\cup W^{-}_{1}(y,u)$. This implies $r\geq |W^{-}_{1}(z_{1},u)\cup W^{-}_{1}(y,u)|$. Since $uy$ and $uz_{1}$ are edges of $F-E(F_{0})$, Claim~\ref{cl:stars} implies they cannot be in the same $1$-factor for any edge of $F-E(F_{0})$ between $w$ and  $W^{-}_{1}(z_{1},u)\cap W^{-}_{1}(y,u)$. Thus, $r\geq |W^{-}_{1}(z_{1},u)\cap W^{-}_{1}(y,u)|+|\{z_{1},y\}|$. This completes the proof of item (\ref{it1}). The proofs for items (\ref{it2}), (\ref{it3}), and (\ref{it4}) have almost identical arguments to this one so we omit their proofs.

If $deg_{G}(z_{1})\leq deg_{G}(u)$, then Claim~\ref{cl:join} says $z_{1}$ is adjacent in $F-E(F_{0})$ to every vertex in $W^{+}_{i}(u,z_{1})$ and $u$ is adjacent in $F-E(F_{0})$ to every vertex in  $W^{-}_{1}(z_{1},u)$. Thus, Claim~\ref{cl:starsObservation} says  $r\geq |W^{+}_{i}(u,z_{1})|+|W^{-}_{1}(z_{1},u)|-1$. We omit the proof of item (\ref{it10}) since it is similar to this one.
\end{proof}
\end{claim}

We are now ready to tackle the heart of the proof. Let us first show that $t\neq 2$.

\begin{claim} $t\neq 2$
\begin{proof}
Suppose $t=2$. In this case, for all $x\in L_{1}\cup L_{2}$, (\ref{eq:forbidden}) implies $e_{F_{0}}(z_{1},A)=e_{F_{0}}(x,A)=0$, and therefore, by our choice of $z_{1}$, $deg_{G}(x)\leq deg_{G}(z_{1})$. Since $e_{F_{0}}(D_{2},S)=0$, $D_{2}$ is $k'$-regular, and therefore, $k'$ is even since $|D_{2}|$ is odd. Thus, we need to show that $r\geq \frac{k'}{2}+3$. Since $W^{-}_{2}(x,z_{1})\neq \emptyset$ and $W^{+}_{1}(z_{1},x)\neq \emptyset$ for all $x\in L_{2}$, Claim~\ref{cl:forbidden} says $r\geq 4$. Thus, we have proved this claim for $k'=2$. From here we will assume $k'\geq 4$ and $r\leq \frac{k'}{2}+2$.

Let $J$ be the fewest edges in $D_{i}$ such that $D_{i}-J$ is bipartite with bi-partitions $X$ and $Y$. Since there are odd number of vertices in $D_{i}$, we may without loss of generality assume $|X|\geq |Y|+1$. Since $J$ is the smallest number of such edges, we see that $J\cap E(X,Y)=\emptyset$. Letting $J_{X}=J\cap E(D_{i}[X])$ and $J_{Y}=J\cap E(D_{i}[Y])$, we see that $k'|X|-2|J_{X}|=e_{D_{i}-J}(X,Y)=k'|Y|-2|J_{Y}|$. This implies $|J|\geq |J_{Y}|-|J_{X}|=\frac{k'}{2}(|X|-|Y|)\geq \frac{k'}{2}$. 

For $x\in L_{2}$, if $deg_{G}(z_{1})=deg_{G}(x)$, $W_{2}^{-}(x,z_{1})=N_{D_{2}}[x]$, or $W_{1}^{+}(z_{1},x)=N_{D_{1}}[z_{1}]$, then Claim~\ref{cl:join} says either $z_{1}$ is adjacent in $F-E(F_{0})$ to the $k'+1$ vertices in $N_{D_{2}}[x]$ or $x$ is adjacent in $F-E(F_{0})$ to the $k'+1$ vertices in $N_{D_{1}}[z_{1}]$. Therefore, $r\geq k'+1=\frac{k'}{2}+\frac{k'}{2}+1\geq \frac{k'}{2}+3$. Thus, we may assume for all $x\in L_{2}$ that $deg_{G}(z_{1})>deg_{G}(x)$ and both $W_{2}^{+}(x,z_{1})$ and $W_{1}^{-}(z_{1},x)$ are not empty. 

We choose a $z_{2}\in L_{2}$ such that  $|W^{-}_{2}(z_{2},z_{1})|$ is maximized. Since $deg_{G}(z_{2})<deg_{G}(z_{1})$, $z_{2}\notin W^{+}_{2}(z_{2},z_{1})$ and $z_{1}\notin W^{-}_{1}(z_{1},z_{2})$.

For every $x\in L_{2}$, $r\geq |W^{-}_{2}(x,z_{1})|+|W^{+}_{1}(z_{1},x)|-1$ by Claim~\ref{cl:joinDistinct} item (\ref{it10}). Thus, we may assume $|W^{-}_{2}(z_{2},z_{1})|+|W^{+}_{1}(z_{1},z_{2})|\leq \frac{k'}{2}+3$. Letting $f= |W^{+}_{1}(z_{1},z_{2})|$ we see that $|W^{-}_{2}(z_{2},z_{1})|\leq \frac{k'}{2}+3-f$. Therefore, \[|W^{+}_{2}(z_{2},z_{1})|=|N_{D_{2}}[z_{2}]|- |W^{-}_{2}(z_{2},z_{1})|\geq k'+1- \left(\frac{k'}{2}+3-f\right) = \frac{k'}{2}+f-2\] and \[|W^{-}_{1}(z_{1},z_{2})|=|N_{D_{1}}[z_{1}]|-|W^{+}_{1}(z_{1},z_{2})|=k'+1-f.\] By Claim~\ref{cl:joinDistinct} items (\ref{it1}) and (\ref{it2}), $r\geq |W^{-}_{1}(z_{1},z_{2})|+1$ and $r\geq |W^{+}_{2}(z_{2},z_{1})| +1$. However, $\frac{k'}{2}+2\geq r\geq |W^{+}_{2}(z_{2},z_{1})| +1= \frac{k'}{2}+f-1$ and $\frac{k'}{2}+2\geq r\geq |W^{-}_{1}(z_{1},z_{2})|+1 = k'+2-f$. This implies $2\leq f\leq 3$ and $f\geq \frac{k'}{2}$. Therefore, $4\leq k'\leq 6$.

Suppose $|W^{-}_{2}(z_{2},z_{1})|\leq \frac{k'}{2}$. Therefore, $|E(W^{-}_{2}(z_{2},z_{1}),\{z_{2}\})|<\frac{k'}{2}$ since $z_{2}\in W^{-}_{2}(z_{2},z_{1})$. Previously, we showed that every odd cycle in $D_{2}$ must intersect a set of at least $k'/2$ edges. Thus, there must be an  $x\in L_{2}-\{z_{2}\}$. By our choice of $z_{2}$, $k'/2\geq |W^{-}_{2}(z_{2},z_{1})|\geq |W^{-}_{2}(x,z_{1})|$. Thus, \[\frac{k'}{2} + 1\leq |W^{+}_{2}(z_{2},z_{1})|\leq |W^{+}_{2}(x,z_{1})|\leq \frac{k'}{2}+f-2.\] This implies $f=3$, and therefore, $|W^{-}_{1}(z_{1},z_{2})|\geq k'-2$.
We choose a $w\in N_{D_{1}}(z_{1})$ with the smallest degree in $G$. Since $W_{1}^{-}(z_{1},z_{2})$ and $W_{1}^{-}(z_{1},x)$ are not empty, we may conclude $w\in W_{1}^{-}(z_{1},z_{2})\cap W_{1}^{-}(z_{1},x)$. Therefore, by Claim~\ref{cl:joinDistinct} item (\ref{it4}), \[\frac{k'}{2}+2\geq r\geq \max\{|W_{2}^{+}(x,z_{1})\cup W_{2}^{+}(z_{2},z_{1})|,|W^{+}_{2}(x,z_{1})\cap W^{+}_{2}(z_{2},z_{1})|+|\{z_{2},x\}|\}.\]
This implies $|W_{2}^{+}(x,z_{1})\cap W^{+}_{2}(z_{2},z_{1})|\leq k'/2$. Moreover, since $\frac{k'}{2}+1\leq |W_{2}^{+}(z_{2},z_{1})|\leq |W^{+}_{2}(x,z_{1})|$ and \[|W_{2}^{+}(x,z_{1})\cup W_{2}^{+}(z_{2},z_{1})|=|W_{2}^{+}(x,z_{1})|+ |W^{+}_{2}(z_{2},z_{1})|-|W^{+}_{2}(x,z_{1})\cap W^{+}_{2}(z_{2},z_{1})|,\] we see that \[\frac{k'}{2}+2\geq|W_{2}^{+}(x,z_{1})\cup W^{+}_{2}(z_{2},z_{1})|\geq 2\left(\frac{k'}{2}+1\right)-\frac{k'}{2}=\frac{k'}{2}+2.\] Thus, $|W^{+}_{2}(x,z_{1})- W^{+}_{2}(z_{2},z_{1})|\leq 1$. We now need to do some edge counting. By Claim~\ref{cl:join}, there are $(\frac{k'}{2}+2)|W^{-}_{1}(z_{1},z_{2})|+2$ edges of $F-E(F_{0})$ between $W^{+}_{1}(x,z_{1})\cup W^{+}_{2}(z_{2},z_{1})\cup \{z_{2},x\}$ and $W^{-}_{1}(z_{1},z_{2})\cup \{z_{1}\}$. Since $r\leq \frac{k'}{2}+2$, we may conclude by the pigeonhole principle that there is some $H_{s}$ for $s\geq 3$ that has $|W^{-}_{1}(z_{1},z_{2})|+1$ edges between $W^{+}_{2}(x,z_{1})\cup W^{+}_{2}(z_{2},z_{1})\cup \{z_{2},x\}$ and $W^{-}_{1}(z_{1},z_{2})\cup \{z_{1}\}$. Thus, $z_{1}$ is adjacent in $H_{s}$ to some $x'\in \{x, z_{2}\}$ and every vertex in $W^{-}_{1}(z_{1},z_{2})$ is adjacent in $H_{2}$ to some vertex in $W^{+}_{2}(x,z_{1})\cup W^{+}_{2}(z_{2},z_{1})$. Since $|W^{-}_{1}(z_{1},z_{2})|\geq k'-2\geq 2$ and $|W^{+}_{2}(x,z_{1})- W^{+}_{2}(z_{2},z_{1})|\leq 1$ we may conclude there is a $w\in W^{-}_{1}(z_{1},z_{2})$ and a $w'\in W^{+}_{2}(x',z_{1})$ such that $ww'\in E(H_{s})$. However, the edges $ww'$ and $z_{1}x'$ contradict Claim~\ref{cl:stars}. Thus, $|W^{-}_{2}(z_{2},z_{1})|\geq \frac{k'}{2}+1$, and therefore, $f=2$. 

Since $f=2$ and $z_{1}\in W^{+}_{1}(z_{1},z_{2})$, $k'=4$ and $r=4$. Therefore, $|E(W^{+}_{1}(z_{1},z_{2}),\{z_{1}\})|=1$. Since we showed that every odd cycle in $D_{1}$ must intersects a set of at least two edges, there must be a $y\in L_{1}-\{z_{1}\}$. By our choice of $z_{1}$, $deg_{G}(y)\leq deg_{G}(z_{1})$. Since, $W^{+}_{2}(z_{2},z_{1})\neq \emptyset$, $W^{+}_{2}(z_{2},z_{1})\cap W^{+}_{2}(z_{2},y)\neq \emptyset$. Therefore, Claim~\ref{cl:joinDistinct} item (\ref{it1}) implies \[r\geq \max\{|W^{-}_{1}(z_{1},z_{2})\cup W^{-}_{1}(y,z_{2})|, |W^{-}_{1}(z_{1},z_{2})\cap W^{-}_{1}(y,z_{2})|+|\{z_{1},y\}|\}.\] Since $|W^{-}_{1}(z_{1},z_{2})|=3$ and $r=4$, we may conclude that $|W^{-}_{1}(z_{1},z_{2})\cap W^{-}_{1}(y,z_{2})|\leq 2$. If  $|W^{-}_{1}(y,z_{2})|\geq 3$, then 
$|W^{-}_{1}(z_{1},z_{2})\cup W^{-}_{1}(y,z_{2})|=4$. By Claim~\ref{cl:join}, there are $4|W^{+}_{2}(z_{2},z_{1})|+2$ edges of $F-E(F_{0})$ between $W^{-}_{1}(z_{1},z_{2})\cup W^{-}_{1}(y,z_{2})\cup \{z_{1},y\}$ and $W^{+}_{2}(z_{2},z_{1})\cup \{z_{2}\}$. Since $r=4$, we may conclude by the pigeonhole principle that there is some $H_{s}$ for $s\geq 3$ that has $|W^{+}_{2}(z_{2},z_{1})|+1$ edges between $W^{-}_{1}(z_{1},z_{2})\cup W^{-}_{1}(y,z_{2})\cup \{z_{1},y\}$ and $W^{+}_{2}(z_{2},z_{1})\cup \{z_{2}\}$. Thus, $z_{2}$ is adjacent in $H_{s}$ to some $y'\in \{y, z_{1}\}$, and every vertex in $W^{+}_{2}(z_{2},z_{1})$ is adjacent in $H_{2}$ to some vertex in $W^{-}_{1}(z_{1},z_{2})\cup W^{-}_{1}(y,z_{2})$. Since $|W^{+}_{2}(z_{2},z_{1})|=2$ and $|W^{-}_{1}(z_{1},z_{2})- W^{-}_{1}(y,z_{2})|=1$ we may conclude there is a $w\in W^{+}_{2}(z_{2},z_{1})$ and a $w'\in W^{-}_{1}(y',z_{2})$ such that $ww'\in E(H_{s})$. However, the edges $ww'$ and $z_{2}y'$ contradict Claim~\ref{cl:stars}. Thus, $|W^{-}_{1}(y,z_{2})|\leq 2$, and therefore, $|W^{+}_{1}(y,z_{2})|\geq 2$. Since $z_{2}$ is adjacent in $F-E(F_{0})$ to every vertex in $W_{1}^{+}(z_{1},z_{2})\cup W_{1}^{+}(y,z_{2})$ we may conclude that $|W_{1}^{+}(z_{1},z_{2})\cup W_{1}^{+}(y,z_{2})|\leq 4$, and therefore, there is a $w\in W_{1}^{+}(z_{1},z_{2})\cap W_{1}^{+}(y,z_{2})$. However, Claim~\ref{cl:joinDistinct} item (\ref{it10}) implies $r\geq |W^{-}_{2}(z_{2},z_{1})|+|\{z_{1},y\}|\geq 5$. Thus, $t\neq 2$.\qedhere
\end{proof}
\end{claim}

From this point on we assume $t\geq 3$. We are left with two main cases.

First, we suppose $W^{-}_{1}(z_{1},z_{i})\neq N_{D_{1}}[z_{1}]$ nor $W^{+}_{1}(z_{1},z_{i})\neq N_{D_{1}}[z_{1}]$ for all $z_{i}\in L_{i}$.

For each $i\geq 2$, we choose a $z_{i}\in L_{i}$. We let $z^{-}_{1}\in N_{D_{1}}(z_{1})$ such that $deg_{G}(z^{-}_{1})\leq deg_{G}(u)$ for all $u\in N_{D_{1}}(z_{1})$. Similarly, we let $z^{+}_{1}\in N_{D_{1}}(z_{1})$ such that $deg_{G}(z^{+}_{1})\geq deg_{G}(u)$ for all $u\in N_{D_{1}}(z_{1})$. Thus, $z^{+}\in W^{+}_{1}(z_{1},z_{i})$ and $z^{-}\in W^{-}_{1}(z_{1},z_{i})$ for all $i\in \{2,\ldots,t\}$. By Claim~\ref{cl:join} every vertex in $N_{D_{i}}[z_{i}]$ is adjacent in $F-E(F_{0})$ to either $z_{1}^{-}$ or $z_{1}^{+}$.

By (\ref{eq:forbidden}) \[k'(t-2)\geq \sum_{i=2}^{t}e_{F_{0}}(D_{i},S)\geq \sum_{i=2}^{t}e_{F_{0}}(z_{i},S).\] This implies \[\sum_{i=2}^{t}|N_{D_{i}}(z_{i})|\geq \sum_{i=2}^{t}(k'-e_{F_{0}}(z_{i},S))\geq k'(t-1)-\sum_{i=2}^{t}e_{F_{0}}(z_{i},S)\geq k'(t-1)-k'(t-2)=k'.\] Thus, $\sum_{i=2}^{t}|N_{D_{i}}[z_{i}]|\geq k'+t-1$. This implies either $z^{+}$ or $z^{-}$ is adjacent in $F-E(F_{0})$ to $\lceil\frac{k'+t-1}{2}\rceil$ vertices in $\cup_{i=2}^{t}N_{D_{i}}[z_{i}]$. By Claim~\ref{cl:stars} none of those edges incident with $z^{+}$ or $z^{-}$ can be in the same $1$-factor as $z_{i}z_{1}$ for all $i\geq 2$. Thus, \[r\geq t-1+\bigg\lceil\frac{k'+t-1}{2}\bigg\rceil\geq 2 + \bigg\lceil\frac{k'+2}{2}\bigg\rceil.\]

To finish the proof, we suppose $W^{-}_{1}(z_{1},z_{i})=N_{D_{1}}[z_{1}]$ or $W^{+}_{1}(z_{1},z_{i})=N_{D_{1}}[z_{1}]$ for some $z_{i}\in L_{i}$. Let $M'$ be a matching that misses both $z_{i}$ and $z_{1}$. By Claim~\ref{cl:NumberofComp} $|\mathcal{D}'(z_{i},M')|\geq |N_{F_{0}}(D'(z_{i},M'))\cap A|+1$. Let $q$ be the number of non-trivial components in $\mathcal{D}'(z_{i},M')$.  If $e_{F_{0}}(z_{1},A)\leq \big\lfloor\frac{k'+1}{2}\big\rfloor$, then 
\begin{align}
    k|N_{F_{0}}(D'(z_{i},M'))\cap A|&\geq k(|\mathcal{D'}(z_{i},M')|-q)+qe_{F_{0}}(z_{1},A)\notag\\
    &\geq k(|N_{F_{0}}(D'(z_{i},M'))\cap A|+1-q)+q\bigg\lfloor\frac{k'+1}{2}\bigg\rfloor.\notag
\end{align} This can only happen if $q\geq 2$. Thus, there is a non-trivial component $D_{j}\in \mathcal{D'}(z_{i},M')$ that is not $D_{1}$. Given a $z_{j}\in L_{j}$, Claim~\ref{cl:forbidden} says $z_{i}$ is adjacent in $F-E(F_{0})$ to at least two vertices in $D_{j}$. This implies \[r\geq |N_{D_{1}}(z_{1})|+2\geq k'-e_{F_{0}}(z_{1},A)+3\geq k'-\bigg\lfloor\frac{k'+1}{2}\bigg\rfloor+3=\bigg\lceil\frac{k'+1}{2}\bigg\rceil+2.\] Thus, we are left with the case $e_{F_{0}}(z_{1},A)\geq \big\lfloor\frac{k'+1}{2}\big\rfloor+1$. 

By our choice of $z_{1}$, $e_{F_{0}}(z_{i},A)\geq e_{F_{0}}(z_{1},A)\geq \big\lfloor\frac{k'+1}{2}\big\rfloor+1$. We stated earlier that $|D'|\geq |S|+1$ and the first $t$ components of $D'$ are non-trivial. Thus, \[k'|S|\geq e_{F_{0}}(D',S)\geq k'(|S|+1-(t-1))+(t-1)\bigg(\bigg\lfloor\frac{k'+1}{2}\bigg\rfloor+1\bigg).\] Solving for $t-1$ we see that $t-1\geq 3$, and therefore, $t\geq 4$.

Since $e_{F_{0}}(D_{1},A)\leq k'-1$, there exists an $x\in L_{i}$ such that \[\bigg\lfloor\frac{k'+1}{2}\bigg\rfloor+1\leq e_{F_{0}}(z_{1},A)\leq e_{F_{0}}(x,A)\leq (k'-1)/|L_{i}|.\] This implies $k'\geq 4$ and $L_{i}=\{z_{i}\}$ for all $i$. Thus, every odd cycle $C_{j}$ in $D_{i}$ passes through $z_{i}$ with $z_{j}z^{+}_{i}$ and $z^{-}_{j}z_{i}$ being edges of $C$ such that $deg_{G}(z_{j}^{-})\leq deg_{G}(z_{i}
)\leq deg_{G}(z_{j}^{+})$. Note that if $z_{i}\in W^{-}_{i}(z_{i},z_{1})$, then $z^{-}_{j}\in W^{-}_{i}(z_{i},z_{1})$. Likewise, if $z_{i}\in W^{+}_{i}(z_{i},z_{1})$, then $z^{+}_{j}\in W^{+}_{i}(z_{i},z_{1})$. Thus, there is a smallest set $J_{i}=\{z_{i}u^{(i)}_{1},\ldots, z_{i}u^{(i)}_{l_{i}}\}$ of $l_{i}$ edges such that $E(C)\cap J_{i}\neq \emptyset$ for every odd cycle $C$ in $D_{i}$ and $\{u^{(i)}_{1},\ldots, u^{(i)}_{l_{i}}\}$ is a subset of $W_{i}^{-}(z_{i},z_{1})$ when $deg_{G}(z_{1})\geq deg_{G}(z_{i})$ or a subset of $W_{i}^{+}(z_{i},z_{1})$ when $deg_{G}(z_{1})\leq deg_{G}(z_{i})$. The key take away here is that $z_{1}$ is adjacent in $F-E(F_{0})$ to every vertex in $\{z_{i},u^{(i)}_{1},\ldots, u^{(i)}_{l_{i}}\}$ for all $i\geq 1$. Thus, we are left with the task of establishing a lower bound on $\sum_{i=2}^{t}l_{i}$.

Since $J_{i}$ is made up of $l_{i}$ distinct edges, we deduce that $l_{i}\leq k'-e_{F_{0}}(z_{i},A)$ for all $i\geq 2$. For each $D_{i} \in \mathcal{D'}$, $D_{i}-J_{i}$ must be a bipartite graph with partitions $X_{i}$ and $Y_{i}$ since every odd cycle of $D_{i}$ shares an edge with $J_{i}$. We may assume without loss of generality that $|X_{i}|\geq |Y_{i}|+1$ since $D_{i}$ has an odd number of vertices. Let $w_{i}$ be the number of edges of $J_{i}$ with both ends in $X_{i}$, and $w'_{i}$ be the number of edges of $J_{i}$ with both ends in $Y_{i}$.  Since $J_{i}$ is the smallest such set, we may deduce that no edge of $J_{i}$ is incident with both $X_{i}$ and $Y_{i}$. Thus, $w'_{i}=|J_{i}|-w_{i}$. Moreover, since every edge of $J_{i}$ is incident with $z_{i}$ we may conclude that either $w_{i}=0$ or $w'_{i}=0$.  Thus, \[k'|Y_{i}|-2w'_{i}-e_{F_{0}}(Y_{i},S)= e_{F_{0}}(X_{i},Y_{i})= k'|X_{i}|-2w_{i}-e_{F_{0}}(X_{i},S).\] Rewriting this we see that \[e_{F_{0}}(X_{i},S)=k'(|X_{i}|-|Y_{i}|)+2(w'_{i}-w_{i})+e_{F_{0}}(Y_{i},S)\geq k'-2w_{i}+e_{F_{0}}(Y_{i},S)\geq k'-2l_{i}+e_{F_{0}}(Y_{i},S) .\] Thus, \[e_{F_{0}}(D_{i},S)=e_{F_{0}}(X_{i},S)+e_{F_{0}}(Y_{i},S)\geq k'-2l_{i}+2e_{F_{0}}(Y_{i},S)\geq k' - 2l_{i}.\] Therefore, \[\sum_{i=2}^{t}e_{F_{0}}(D_{i},S)\geq (t-1)k'-2\sum_{i=2}^{t}l_{i}.\] Combining this last inequality with (\ref{eq:forbidden}) we see that $k'(t-2)\geq (t-1)k'-2\sum_{i=2}^{t}l_{i}$. This implies $\sum_{i=2}^{t}l_{i}\geq k'/2$, and therefore, $z_{1}$ is adjacent in $F-E(F_{0})$ to at least \[\frac{k'}{2}+t-1\geq \frac{k'}{2}+3=\bigg\lceil\frac{k'+1}{2}\bigg\rceil+2\] vertices in $D'$. Thus, $r\geq \lceil\frac{k'+1}{2}\big\rceil+2$. \qedhere
\end{proof}

\providecommand{\bysame}{\leavevmode\hbox to3em{\hrulefill}\thinspace}
\providecommand{\MR}{\relax\ifhmode\unskip\space\fi MR }
\providecommand{\MRhref}[2]{%
  \href{http://www.ams.org/mathscinet-getitem?mr=#1}{#2}
}

\providecommand{\bysame}{\leavevmode\hbox to3em{\hrulefill}\thinspace}
\providecommand{\MR}{\relax\ifhmode\unskip\space\fi MR }
\providecommand{\MRhref}[2]{%
  \href{http://www.ams.org/mathscinet-getitem?mr=#1}{#2}
}
\providecommand{\href}[2]{#2}

\end{document}